\documentclass[11pt,a4paper]{amsart} 
\usepackage[english]{babel}
\usepackage{amssymb,xspace,enumerate,graphics,ifthen}
\usepackage{amstext}
\theoremstyle{plain}
\usepackage{amsbsy,amssymb,amsfonts,latexsym}

\date{\today}
\title{Division of holomorphic functions and growth conditions}
\author{William ALEXANDRE}
\address{Laboratoire Paul Painlev\'e U.M.R. CNRS 8524, U.F.R. de
Math\'ematiques,  cit\'e scientifique, Universit\'e Lille 1, F59 655 Villeneuve d'Ascq Cedex, France.}
\email{ william.alexandre@math.univ-lille1.fr}
\thanks{The first author is partially supported by A.N.R. BL-INTER09-CRARTIN}
\author{Emmanuel MAZZILLI}
\address{Laboratoire Paul Painlev\'e U.M.R. CNRS 8524, U.F.R. de
Math\'ematiques,  cit\'e scientifique, Universit\'e Lille 1, F59 655 Villeneuve d'Ascq Cedex, France.}
\email{ emmanuel.mazzilli@math.univ-lille1.fr}
\date{}
\newtheorem{theorem}{Theorem}[section]
\newtheorem{lemma}[theorem]{Lemma}

\newtheorem{proposition}[theorem]{Proposition}
\newtheorem{corollary}[theorem]{Corollary}

\def \pint {\vbox{ \hbox to 5 pt {\hfil \vrule height 4pt}\hrule}\hskip 3pt}
\def\leqs{\lesssim}
\def\geqs{\gtrsim}
\def\eqs{\eqsim}
\def\cc{\mathbb{C}}
\def\rr{\mathbb{R}}
\def\nn{\mathbb{N}}

\def \qed {\hbox{\hskip 5pt} \vbox{\hrule \hbox to 5pt 
{\vrule height 4.2pt \hfil \vrule}\hrule}}
\def \pint {\vbox{ \hbox to 5 pt {\hfil \vrule height 4pt}\hrule}\hskip 3pt}

\newcommand{\cal}{\mathcal}
\newcommand{\wrt}{with respect to }

\newcommand{\diffp}[2]{\frac{\partial #1}{\partial #2}}
\newcommand{\mlabel}[1]{\label {#1}}
\renewcommand{\over}[2]{\genfrac{}{}{0pt}{}{#1}{#2}}

\newcommand{\zjk}{z_{j,k}}
\newcommand{\p}[2]{{\cal P}_{#1}(#2)}

\newcommand{\pk}[2]{{\cal P}_{{\ifthenelse{#1=1}{ }{#1}}\kappa|\rho({#2})|}({#2})}
\newcommand{\rzjk}{|\rho(\zjk)|}

\newcommand{\pr} {\noindent{\it Proof:} }
\newcommand{\ko} {Koranyi }
\newcommand{\kb}[2] {{\cal P}_{#1}(#2)}

\newcommand{\oo}{{\cal O}}
\begin{document}
\pagestyle{plain}

\begin{abstract}
Let $D$ be a strictly convex domain of $\cc^n$, $f_1$ and $f_2$ be two holomorphic functions defined on a neighborhood of $\overline D$ and set $X_l=\{z,\ f_l(z)=0\}$, $l=1,2$. Suppose that $X_l\cap bD$ is transverse for $l=1$ and $l=2$, and that $X_1\cap X_2$ is a complete intersection. We give necessary conditions when $n\geq 2$ and sufficient conditions when $n=2$ under which a function $g$ to be written as $g=g_1f_1+g_2f_2$ with $g_1$ and $g_2$ in $L^q(D)$, $q\in [1,+\infty)$, or $g_1$ and $g_2$ in $BMO(D)$. In order to prove the sufficient condition, we explicitly write down the functions $g_1$ and $g_2$ using integral representation formulas and new residue currents.
\end{abstract}
\maketitle
\section{Introduction}
In this article, we are interested in ideals of holomorphic functions and corona type problems. More precisely, if $D$ is a domain of $\cc^n$ and $f_1,\ldots, f_k$ are $k$ holomorphic functions defined in a neighborhood of $\overline D$, we are looking for condition(s), as close as possible to being necessary and sufficient, under which a function $g$, holomorphic on $D$, can be written as 
\begin{eqnarray}
 g&=&f_1g_1+\ldots+ f_kg_k,\label{eq0}
\end{eqnarray}
with $g_1,\ldots, g_k$ holomorphic on $D$ and satisfying growth conditions at the boundary of $D$. This kind of problem has been widely studied by many authors under different assumptions.
\par\medskip
When $D$ is strictly pseudoconvex and when $f_1,\ldots, f_k$ are holomorphic and bounded functions on $D$, which satisfy  $|f|^2=|f_1|^2+\ldots+|f_k|^2\geq\delta^2 >0$, for a given holomorphic and bounded function $g$, finding functions $g_1,\ldots, g_k$ bounded on $D$  is a question known as the Corona Problem. When $D$ is the unit ball of $\cc$, the Corona Problem was solved in 1962 by Carleson in \cite{Car}. This question is still open for $n>1$, even for two generators $f_1$ and $f_2$, and even when $D$ is the unit ball of $\cc^n$.
\par
For $p\in [1,+\infty)$, we denote by $H^p(D)$ the Hardy space of $D$.
When $n>1$, $k=2$ and $|f|\geq \delta >0$, Amar proved in \cite{Ama} that for any $g\in H^p(D)$, (\ref{eq0}) can be solved with $g_1$ and $g_2$ in $H^p(D)$. Andersson and Carlsson in \cite{AC1} generalized this result to any strictly pseudoconvex domain in $\cc^n$ and to any $k\geq 2$ and also obtained the $BMO$-result already announced by Varopoulos in \cite{Var}. In \cite{AC3}, they studied the dependence of the $g_i$'s on the lower bound $\delta$ of $|f|$ and they explicitly obtained a constant  $c_\delta $ such that for all $i$, $\|g_i\|_{H^p(D)}\leq c_\delta \|g\|_{H^p(D)}$. Of course $c_\delta $ goes to infinity when $\delta$ goes to $0$.  In \cite{AB}, when $|f|$ does not have a positive lower bound, Amar and Bruna formulated a sufficient condition in term of the admissible maximum function of $|f|^{−2}\left|\log|f|\right|^{2+\varepsilon }$, $\varepsilon >0,$ under which the $g_i$'s belong to $H^p(D)$.

The corona problem was also studied in the case of the Bergman space $A^p(D)$, the space of holomorphic functions which belong to $L^p(D)$, and in the case of the Zygmund space $\Lambda_\gamma(D)$ by Krantz and Li in \cite{KL}, and in the case of Hardy-Sobolev spaces by Ortega and F\`abrega in \cite{OF}.

In the  above papers, the first step of the proof in the case of two generators $f_1$ and $f_2$, is to find  two smooth functions on $D$, $\varphi _1$ and $\varphi_2$, such that 
\begin{eqnarray}
\varphi _1f_1+\varphi_2 f_2&=&1; \label{eq00}
\end{eqnarray}
and then to solve the equation
\begin{eqnarray}
 \overline\partial \varphi &=&\frac{\overline{f_1}\ \overline\partial\varphi_2-\overline{f_2}\ \overline \partial\varphi _1}{|f_1|^2+|f_2|^2}.\label{eq000}
\end{eqnarray}
Then setting $g_1=g\varphi_1+\varphi f_2$ and $g_2=g\varphi _2-\varphi f_1$, (\ref{eq0}) holds and, provided $\varphi $ belongs to the appropriate space, $g_1$ and $g_2$ will belong to $H^p(D)$, $A^p(D),\ldots$ So the problem is reduced to solve the Bezout equation (\ref{eq00}) and then to solve the $\overline\partial$-equation (\ref{eq000}) with an appropriate regularity.\par
In \cite{AC2}, Andersson and Carlsson used an alternative technique. They constructed a division formula $g=f_1 T_1(g)+\ldots+f_k T_k(g)$ where for all $i$, $T_i$ is a well chosen Berndtsson-Andersson integral operator, and, still under the assumption $|f|\geq \delta >0$, they proved that $T_i(g)$ belongs to $H^p(D)$ (resp. $BMO(D)$) when $g$ belongs to $H^p(D)$ (resp. $BMO(D)$). The same kind of technics was also used in \cite{BCZ} by Bonneau, Cumenge and Z\'eriahi who studied the equation (\ref{eq0}) in Lipschitz spaces and in the space $B_M(D)=\{ g,\ \|g\|_{B_M(D)}=\sup_{z\in D} \left(|g(z)| d(z,bD)^M\right)<\infty\}$. In this later work, the generators $f_1,\ldots, f_k$ may have common zeroes but $\partial f_1\wedge \ldots\wedge\partial f_k$ can not vanish on $bD\cap \{z,\ f_1(z)=\ldots=f_k(z)=0\}$. 

The case of generators having common zeroes has also been investigated by Skoda in \cite{Sko} for weighted $L^2$-spaces. Using and adapting the $L^2$-techniques developed by H\"ormander, for $D$ pseudoconvex in $\cc^n$, $\psi$ a plurisubharmonic weight on $D$,  $f_1,\ldots, f_k$ holomorphic in $D$, $q=\inf(n,k)$, $\alpha>1$ and $g$ holomorphic in $D$ such that $\int_D\frac{|g|^2}{|f|^{2\alpha q+2}} e^{-\psi}<\infty$, Skoda showed that there exist $g_1,\ldots, g_k \in {\cal O}(D)$ such that (\ref{eq0}) holds and such that for all $i$, $\int_D \frac{|g_i|^2}{|f|^{2\alpha q}}e^{-\psi}\leq \frac{\alpha}{\alpha-1}\int_D\frac{|g|^2}{|f|^{2\alpha q+2}} e^{-\psi}$. Moreover the result also holds when $k$ is infinite and there is no restriction on $\partial f_1,\ldots, \partial f_k$. However, if one take $g=f_1$ for example, $g$ does not satisfy the assumption of Skoda's theorem in general.
\par\medskip
In this article we restrict ourself to a strictly convex domain $D$ of $\cc^n$ and we consider the case of two generators $f_1$ and $f_2$, holomorphic in a neighborhood of $\overline D$. We denote by $X_1$ the set $X_1=\{z,\ f_1(z)=0\}$, and by $X_2$ the set $X_2=\{z,\ f_2(z)=0\}$. We assume that the intersections $X_1\cap bD$ and  $X_2\cap bD$ are transverse in the sense of tangent cones and that $X_1\cap X_2$ is a complete intersection. We seek assumptions on $g$, holomorphic in $D$, as close as possible to being necessary and sufficient, under which we can write $g$ as $g=g_1f_1+g_2f_2$ with $g_1$ and $g_2$ holomorphic and belonging to $BMO(D)$ or $L^q(D)$, $q\in [1,+\infty)$. 

Let us write $D$ as $D=\{z\in\cc^n,\ \rho(z)<0\}$ where $\rho$ is a smooth strictly convex function defined on $\cc^n$ such that the gradient of $\rho$ does not vanish in a neighborhood $\cal U$ of $bD$. We denote by $D_r$, $r\in\rr$, the set $D_r=\{z\in\cc^n,\ \rho(z)<r\}$, by $\eta_\zeta$  the outer unit normal to $bD_{\rho(\zeta)}$ at a point $\zeta\in{\cal U}$ and by $v_\zeta$ a smooth unitary complex vector field tangent at $\zeta$ to $bD_{\rho(\zeta)}$.
As a first result, we show~:
\begin{theorem}\label{main_result} Let $D$ be a strictly convex domain of $\cc^2$, 
$f_1$ and $f_2$ be two holomorphic functions defined on a neighborhood of $\overline D$ and set $X_l=\{z,\ f_l(z)=0\}$, $l=1,2$. Suppose that $X_l\cap bD$ is transverse for $l=1$ and $l=2$, and that $X_1\cap X_2$ is a complete intersection.\\
Then there exist two integers $k_1,k_2\geq 1$ depending only from $f_1$ and $f_2$ such that if $g$ is any holomorphic function on $D$ which belongs to the ideal generated by $f_1$ and $f_2$ and for which there exist two $C^\infty$ smooth functions $\tilde g_1$ and $\tilde g_2$ such that
\begin{enumerate}[(i)]
 \item\label{mth1} $g=\tilde g_1f_1+\tilde g_2f_2$ on $D$,
 \item\label{mth2} there exists $N\in\nn$ such that $|\rho|^N\tilde g_1$ and $|\rho|^N\tilde g_2$ vanish to order $k_2$ on $bD$,
 \item\label{mth3}there exists $q\in[1,+\infty]$ such that for $l=1,2$, $\left|\diffp{^{\alpha+\beta} \tilde g_l}{\overline{\eta_\zeta}^\alpha\partial\overline{v_\zeta}^\beta}\right| |\rho|^{\alpha+\frac\beta2}$ belongs to $L^q(D)$ for all non-negative integers $\alpha$ and $\beta$ with $\alpha+\beta\leq k_1$,
\end{enumerate}
then there exist two holomorphic functions $g_1,g_2$ on $D$ which belong to $L^q(D)$ if $q<+\infty$ and to $BMO(D)$ if $q=+\infty$, such that $g_1f_1+g_2f_2=g$ on $D$.
\end{theorem}
The number $k_1$ and $k_2$ are almost equal to the maximal order of the singularity of $X_1$ and $X_2$. The functions $g_1$ and $g_2$ will be obtained via integral operators acting on $\tilde g_1$ and $\tilde g_2$. These operators are a combination of a Berndtsson-Andersson kernel and of two (2,2)-currents $T_1$ and $T_2$ such that $f_1T_1+f_2T_2=1$. So instead of first solving the Bezout equation (\ref{eq00}) in the sense of smooth functions, we solve it in the sense of currents and then, instead of solving a $\overline\partial$-equation, we ``holomorphy'' the smooth solutions $\tilde g_1$ and $\tilde g_2$ of the equation $g=\tilde g_1f_1+\tilde g_2f_2$ with integral operators using $T_1$ and $T_2$. These operators can be constructed starting from  any currents $\tilde T_1$ and $\tilde T_2$ such that $f_1\tilde T_1+f_2\tilde T_2=1$ (see section \ref{section4}). However, not all such currents will give operators such that $g_1$ and $g_2$ belongs to $L^q(D)$ or $BMO(D)$; they have to be constructed taking into account the behavior of $f_1$ and $f_2$ and more precisely the interplay between $X_1$ and $X_2$ (see section \ref{construction}). Moreover, if $\tilde g_1$ and $\tilde g_2$ are already holomorphic and satisfy the assumptions $(i)-(iii)$ of Theorem \ref{main_result}, then $g_1=\tilde g_1$ and $g_2=\tilde g_2$.

Observe that in Theorem \ref{main_result}, we do not make any assumption on $f_1$ or $f_2$ excepted that the intersection $X_1\cap bD$ and $X_2\cap bD$ are transverse in the sense of tangent cones, and that $X_1\cap X_2$ is a complete intersection. This later assumption can be removed provided we add a fourth assumption on $\tilde g_1$ and $\tilde g_2$. If we moreover assume that\\[5pt]
{\it (iv)} $\diffp{^{\alpha+\beta} \tilde g_1}{\overline{\eta_\zeta}^\alpha\partial\overline{v_\zeta}^\beta}=0$ on $X_2\cap D$ and $\diffp{^{\alpha+\beta} \tilde g_2}{\overline{\eta_\zeta}^\alpha\partial\overline{v_\zeta}^\beta}=0$ on $X_1\cap D$ for all non negative integers $\alpha$ and $\beta$ with $\alpha+\beta\leq k_1$,\\[5pt]
then Theorem \ref{main_result} also holds whenever $X_1\cap X_2$ is not complete. However, it then becomes very difficult to find $\tilde g_1$ and $\tilde g_2$ which satisfy this fourth assumption, excepted if $X_1\cap X_2$ is actually complete.

Indeed, the main difficulty in order to be able to apply Theorem \ref{main_result} is to find the two functions $\tilde g_1$ and $\tilde g_2$ satisfying {\it (\ref{mth1})-(\ref{mth3})}. The canonical choice when $|f|\geq\delta >0$ is to set $\tilde g_1=g\overline {f_1}  |f|^{-2}$ and $\tilde g_2=g\overline {f_2}  |f|^{-2}$. If $|f|\geq\delta >0$ and if $g$ belongs to $L^q(D)$, then $\tilde g_1$ and $\tilde g_2$ will satisfy {\it (\ref{mth1})-(\ref{mth3})} and we can then apply Theorem \ref{main_result}. However, if $|f|$ does not admit a positive lower bound this will not be necessarily the case. For example, when $D=\{z\in\cc^2,\ |z_1-1|^2+|z_2|^2<1\}$, $f_1(z)=z_2$, $f_2(z)=z_2-z^2_1$ and $g=f_1$, we can obviously find $\tilde g_1$ and $\tilde g_2$ which satisfy the assumption of Theorem \ref{main_result} but if we make the canonical choices for $\tilde g_1$ and $\tilde g_2$, they do not fulfill {\it (\ref{mth3})} for $q=\infty$.

Therefore the question of the existence of $\tilde g_1$ and $\tilde g_2$ may itself become a problem that we have to solve.  Using first Koranyi balls, we will reduce this global question to a local one and then, using divided differences, 
we will give numerical conditions under which there indeed exist functions satisfying the hypothesis of Theorem \ref{main_result}. We will also prove that these conditions are necessary in order to solve (\ref{eq0}) with the $g_i$'s belonging to $L^q(D)$, $q\in [1,+\infty]$, even in $\cc^n$. This 
leads us to an effective way of construction of the solutions of (\ref{eq0}) belonging to $L^q(D)$ or $BMO(D)$. 
\par\medskip
The Koranyi balls are defined as follows. We call the coordinates
system centered at $\zeta$ of basis $\eta_{\zeta}, v_{\zeta}$ the \ko coordinates
 at $\zeta$. We denote by $(z_1^*,z_2^*)$ the coordinates of a
point $z$ in the \ko coordinates at $\zeta$. The \ko ball centered
in $\zeta$ of radius $r$ is the set $\kb r {\zeta}:=\{\zeta+\lambda\eta_{\zeta}+\mu
v_{\zeta},\ |\lambda|<r,\ |\mu|<r^{\frac12}\}$.
The following theorem enables us to go from a local division formula in $L^\infty$ to a global division formula in $BMO$.
\begin{theorem}\label{th3}
Let $D$ be a strictly convex domain of $\cc^2$, 
 $f_1$ and $f_2$ be two holomorphic functions defined on a neighborhood of $\overline D$ and set $X_l=\{z,\ f_l(z)=0\}$, $l=1,2$. Suppose that $X_1\cap bD$ and $X_2\cap bD$ are transverse, and that $X_1\cap X_2$ is a complete intersection.\\
Let $g$ be a function holomorphic on $D$ and assume that there exists $\kappa>0$ such that for all $z\in D$, there exist two functions $\hat g_1$ and $\hat g_2$, depending on $z$, $C^\infty$-smooth on ${\cal P}_{\kappa |\rho(z)|}(z)$, such that
 \begin{enumerate}[(a)]
  \item \label{th3i} $g=\hat g_1f_1+\hat g_2 f_2$ on ${\cal P}_{\kappa |\rho(z)|}(z)$;
  \item \label{th3ii} for all non negative integers $\alpha$, $\beta$, $\overline\alpha$ and $\overline\beta$, there exist $c>0$, not depending on $z$, such that $\sup_{{\cal P}_{\kappa |\rho(z)|}(z)}\left|\diffp{^{\alpha+\overline\alpha+\beta+\overline\beta}\hat g_l}{{z^*_1}^\alpha\partial{z_2^*}^\beta\partial \overline{z^*_1}^{\overline\alpha}\partial \overline{\zeta_2^*}^{\overline\beta}}\right|\leq c$ for $l=1$ and $l=2$. 
 \end{enumerate}
Then there exist two smooth functions  $\tilde g_1$ and $\tilde g_2$  which satisfy the assumptions {\it (\ref{mth1})-(\ref{mth3})} of Theorem \ref{main_result} for $q=+\infty$.
\end{theorem}
An analogous theorem holds true in the $L^q$-case (see Theorem \ref{th4}). We observe that if, for all $z\in D$, there exist two functions $\hat g_1$ and $\hat g_2$, holomorphic and bounded on $\pk2z$ by a constant $c$ which does not depend from $z$, and such that $g=\hat g_1f_1+\hat g_2f_2$ on $\pk2z$, then Cauchy's inequalities implies that $\hat g_1$ and $\hat g_2$ satisfy the assumption of Theorem \ref{th3} on $\pk1z$ for all $z$. Therefore  Theorem \ref{th3} implies that the global solvability of (\ref{eq0}) in the $BMO$ space of $D$ is nearly equivalent to its uniform local solvability. In order to prove Theorem \ref{th3}, we will cover $D$ with \ko balls and using a suitable partition of unity, we will glue together the $\hat g_1$ and $\hat g_2$ which we got on each ball. We point out that when we glue together the local $\hat g_1$'s, excepted if $X_1\cap X_2$ is a complete intersection, in general the ``fourth'' assumption {\it (iv)} of Theorem \ref{main_result} is not satisfied. This is why we chose to present Theorem \ref{main_result} as we did.
\par\medskip
When looking for necessary conditions in order to solve (\ref{eq0}) with $g_1$ and $g_2$ bounded, we first observe that $g$ is  trivially bounded by $\max(\|g_1\|_{L^\infty},\|g_2\|_{L^\infty}) (|f_1|+|f_2|)$. Therefore, in order for $g$ to be written as  $g=g_1f_1+g_2f_2$ with $g_1$ and $g_2$ bounded, it is necessary that $\frac{|g|}{|f_1|+|f_2|}$ be bounded. However this condition alone is not sufficient in general. Consider for example the ball $D:=\{z\in\cc^2,\ \rho(z)=|z_1-1|^2+|z_2|^2-1<0\}$, $f_1(z)=z^2_2$, $f_2(z)=z^2_2-z^q_1$ and $g(z)=z_1^{\frac q2}z_2$ where $q\geq 3$ is an odd integer. Then $g(z)={z_2}{z_1^{-\frac q2}} f_1(z)-{z_2}{z_1^{-\frac q2}}f_2(z) $, so $g$ belongs to the ideal generated by $f_1$ and $f_2$, and $\frac{|g|}{|f_1|+|f_2|}$ is bounded on $D$ by $\frac32$; in particular, the classical choice $\tilde g_1=\frac{g \overline {f_1}}{|f_1|^2+|f_2|^2}$ and $\tilde g_2=\frac{g \overline {f_2}}{|f_1|^2+|f_2|^2}$ are smooth and bounded on $D$. However, (\ref{eq0}) can not be solved with $g_1$ and $g_2$ bounded on $D$. In order to see this, a good tool is divided differences. Indeed, on the one hand, if $g=g_1f_1+g_2f_2$, then  $g_1=g\cdot f^{-1}_1$ on $X_2\setminus X_1$. On the other hand, if $g_1$ is bounded, 
for all $z\in D$, all unit vector $v$ tangent to $bD_{-\rho (z)}$ at $z$, all complex numbers $\lambda _1$ and $\lambda _2$ with $\rho(z+\lambda _1v)<\frac12{\rho(z)}$ and $\rho(z+\lambda _2v)<\frac12{\rho(z)}$, the divided difference $\frac{g_1(z+\lambda_1v)-g_1(z+\lambda_2v)}{\lambda _1-\lambda _2}$ behaves like the derivative $\diffp{g_1}{v}$ at some point $z+\mu v$ where $\mu $ is an element of the segment $[\lambda _1,\lambda _2]$ (see \cite{Mon}). 
Cauchy's inequalities then imply that, up to a uniform multiplicative constant, $\frac{g_1(z+\lambda_1v)-g_1(z+\lambda_2v)}{\lambda _1-\lambda _2}$ is bounded by $\|g_1\|_{L^\infty(D)} |\rho(z)|^{-\frac12}$.

So when we compute the divided differences of $g_1$ at points $z+\lambda _1v$ and $z+\lambda _2v$ which belong to  $X_2\setminus X_1$, whatever $g_1$ and $g_2$ may be,  we actually compute the divided difference of $g\cdot f^{-1}_1$; if $g_1$ is bounded, this divided difference  times $|\rho(z)|^{\frac12}$ must be bounded by some uniform constant. But in our example, this is not the case
because for small $\varepsilon>0$, setting $z=(\varepsilon ,0)$, $v=(0,1)$, $\lambda_1=\varepsilon ^{\frac q2}$ and $\lambda_2=-\varepsilon ^{\frac q2}$, we have that 
$ \frac{(g\cdot f_1^{-1}) (z+\lambda_1v)-(g\cdot f_1^{-1}) (z+\lambda_2v)}{\lambda _1-\lambda _2} |\rho(z)|^{\frac12}=\varepsilon ^{\frac{1-q}2}$ which is unbounded when $\varepsilon $ goes to zero.

In $\cc^n$, we will prove that the divided differences of any order of $g\cdot{f_1}^{-1}$ and $g\cdot{f_2}^{-1}$  must satisfy some boundedness properties when (\ref{eq0}) is solvable with $g_1$ and $g_2$ in $L^q(D)$, $q\in [1,+\infty]$ (see Theorems \ref{th1} and \ref{th5} for precise statements). Conversely, in $\cc^2$, if those boundedness properties are satisfied, up to an error term we will be able to construct by interpolation $\hat g_1$ and $\hat g_2$ on any \ko balls which satisfy the assumptions of Theorem \ref{th3}; applying Theorem \ref{th3}, we will then prove that there exist two functions $g_1$ and $g_2$ holomorphic on $D$, belonging to $BMO(D)$ or $L^q(D)$, $q\in [1,+\infty)$, such that $g=g_1f_1+g_2f_2$ (see Theorem \ref{th2} and \ref{th6}).
\par\medskip
The article is organized as follows. In Section \ref{section2}, we recall some tools needed for the construction and the estimation of the division formula. In Section \ref{construction}, we construct the currents which enable us to construct our division formula in Section \ref{section4}. In Section \ref{section5} we establish Theorem \ref{main_result} and finally, in Section \ref{section6}, we prove the theorems related to local division in the $L^\infty$ and $L^q$ case.
       
\section{Notations and tools}\label{section2}
\subsection{Koranyi balls}
The Koranyi balls centered at a point $z$ in $D$ have properties linked with distance from $z$ to the boundary of $D$ in a direction $v$. 
For $z\in \cc^n$, $v$ a unit vector in $\cc^n$, and $\varepsilon>0$, the distance from $z$ to $bD_{\rho(z)+\varepsilon}$ in the direction $v$ 
is defined by 
$$\tau(z,v,\varepsilon)=\sup\{\tau>0,\ \rho(z+\lambda v)-\rho(z)<\varepsilon\text{ for all } \lambda\in\cc,\ |\lambda|<\tau\}.$$
Thus $\tau(z,v,\varepsilon)$ is the maximal radius $r>0$ such that the disc 
$\Delta_{z,v}\left(r\right)=\{z+\lambda v,\ |\lambda|<r\}$ is in $D_{\rho(z)+\varepsilon}$; if $v$ is a tangent vector to $bD_{\rho(z)}$ at $z$, then $\tau(z,v,\varepsilon)$ is comparable to $\varepsilon^{\frac12}$ and $\tau(z,\eta_z,\varepsilon )$ is comparable to $\varepsilon $. \\
Before we recall the properties of the  Koranyi balls we will need, we adopt the following notation. We write $A\leqs B$ if there exists some constant $c>0$ such that $A\leq cB$. Each time we will mention from which parameters $c$ depends. We will write $A\eqs B$ if $A\leqs B$ and $B\leqs A$ both holds.
\begin{proposition}\mlabel{propII.0.1} There exists a neighborhood $\cal U$ of $bD$ and positive real numbers $\kappa$ and $c_1$ such that 
\begin{enumerate}[(i)]
 \item  for all $\zeta\in {\cal U}\cap D$, ${\cal P}_{4\kappa|\rho(\zeta)|}(\zeta)$ is included in $D$.
 \item  for all $\varepsilon>0$, all $\zeta,z\in {\cal U}$, $\p\varepsilon\zeta\cap\p\varepsilon z\neq
\emptyset$ implies $\p\varepsilon z\subset \p{c_1\varepsilon}\zeta$.
\item  for all $\varepsilon>0$ sufficiently small, all
$z\in {\cal U}$, all $\zeta\in \p{\varepsilon}z$ we have $|\rho(z)-\rho(\zeta )|\leq
c_1 \varepsilon$.
\item for all $\varepsilon >0$, all unit vectors $v\in\cc^n$, all $z\in{\cal U}$ and all $\zeta\in{\cal P}_\varepsilon (z)$, $\tau(z,v,\varepsilon )\eqs\tau(\zeta ,v,\varepsilon )$ uniformly with respect to $\varepsilon ,$ $z$ and $\zeta $.
\end{enumerate}
\end{proposition}
For $\cal U$ given by Proposition \ref{propII.0.1} and $z$ and $\zeta$ belonging to $\cal U$, we set $\delta(z,\zeta)=\inf\{\varepsilon>0, \zeta\in \p\varepsilon z\}$. Proposition \ref{propII.0.1} implies  that $\delta$ is a pseudo-distance in the
following sense: 
\begin{proposition}\label{propII.0.2}
For $\cal U$ and $c_1$ given by Proposition \ref{propII.0.1} and for all $z,\ \zeta$ and $\xi$ belonging to $\cal U$  we have
$$\frac1{c_1}\delta(\zeta,z)\leq \delta(z,\zeta)\leq c_1 \delta(\zeta,z)$$
 and
$$\delta(z,\zeta)\leq c_1(\delta(z,\xi)+\delta(\xi,\zeta))$$
\end{proposition}

\subsection{Berndtsson-Andersson reproducing kernel}
\mlabel{secII.0}
Berndtsson-Andersson's kernel will be one of our most important ingredients in the construction of the functions $g_1$ and $g_2$ of Theorem \ref{main_result}. 
We now recall its definition for $D$ a strictly convex domain of $\cc^2$.
We set $h_1(\zeta,z)=-\diffp{\rho}{\zeta_1}(\zeta)$, $h_2(\zeta,z)=-\diffp{\rho}{\zeta_2}(\zeta)$, $h=\sum_{i=1,2} h_id\zeta_i$ and $\tilde{h}=\frac{1}{\rho}
h$. For a $(1,0)$-form $\beta (\zeta,z)=\sum_{i=1,2}\beta_i(\zeta,z)
d\zeta_i$ we set $\langle \beta (\zeta,z),\zeta-z\rangle = \sum_{i=1,2}
\beta_i(\zeta,z)(\zeta_i-z_i)$.
Then we define the Berndtsson-Andersson reproducing kernel by setting for
an arbitrary positive integer $N$, $n=1,2$ and all $\zeta,z\in D$~:
$$P^{N,n}(\zeta,z)=C_{N,n} \left(\frac{1}{1+\langle
\tilde{h}(\zeta,z),\zeta-z\rangle }\right)^{N+n}\left(\overline \partial
\tilde{h}\right)^n,$$
where $C_{N,n}\in \cc$ is a suitable constant. We also set $P^{N,n}(\zeta,z)=0$ for all $z\in D$ and all $\zeta\notin D$. Then the following theorem holds true (see \cite{BA}):
\begin{theorem}\label{BA}
 For all $g\in \oo(D)\cap C^\infty(\overline D)$ we have
$$g(z)=\int_D g(\zeta)P^{N,2}(\zeta,z).$$
\end{theorem}
In order to find an upper bound for this kernel, we will have to write $h$ in the \ko coordinates at some point $\zeta_0$ belonging to $D$. We set   $h_1^*=-\diffp{\rho}{\zeta^*_1}(\zeta)$ and $h_2^*=-\diffp{\rho}{\zeta^*_2}(\zeta)$. Then $h$ is equal to $\sum_{i=1,2}h_i^* d\zeta^*_i$ and satisfies the following Proposition.
\begin{proposition}\mlabel{estiBA}
There exists a neighborhood $\cal U$ of $bD$ such that for all $\zeta\in D\cap {\cal U}$, all $\varepsilon>0$ sufficiently small
and all $z\in \p\varepsilon\zeta$ we have
\begin{enumerate}[(i)]
 \item $|\rho(\zeta)+\langle h(\zeta,z),\zeta-z\rangle|\geqs
\varepsilon+|\rho(\zeta)|+|\rho(z)|$,
 \item $|h^*_1(\zeta,z)|\leqs1$,
\item $|h^*_2(\zeta,z)|\leqs \varepsilon^{\frac{1}{2}}$,
\end{enumerate}
and there exists $c>0$ not depending from $\zeta$ nor from $\varepsilon$ such that
for all $z\in \p\varepsilon\zeta\setminus c\p\varepsilon\zeta$ we have
$$|\langle h(\zeta,z),\zeta-z\rangle|\geqs \varepsilon+|\rho(z)|+|\rho(\zeta)|,$$
uniformly \wrt $\zeta,z$ and $\varepsilon$.
\end{proposition}

\section{Construction of the currents}\label{construction}
In \cite{Maz2}, the following was proved~: If $f_1$ and $f_2$ are two holomorphic functions near the origin  in $\cc^n$, two currents $T$ and $S$  such that $f_1T=1$, $f_2S=\overline\partial T$ and $f_1S=0$ were constructed  on a sufficiently small neighborhood $\cal U$ of $0$. It was also proved that if $T$ and $S$ are any currents satisfying these three hypothesis, then any function $g$ holomorphic on $\cal U$ can be written as $g=f_1g_1+f_2g_2$ on $\cal U$ if and only if $g\overline\partial S=0$. Moreover, $g_1$ and $g_2$ can be explicitly written down using $T$ and $S$.
\par\medskip
Here, when $f_1$ and $f_2$ are holomorphic on a domain $D$, we first want to obtain a decomposition $g=g_1f_1+g_2f_2$ on the whole domain $D$ and then secondly we want to obtain growth estimates on $g_1$ and $g_2$. As a first approach, we could try to globalize the currents $T$ and $S$ of \cite{Maz2} in order to have a global decomposition. However, such an approach would fail to give the growth estimates we want.

In \cite{Maz2}, $f_1$ plays a leading role and $T$ is constructed independently of $f_2$, using only $f_1$. Then $S$ is constructed using $f_1$ and $f_2$. If we assume for example that $f_1$ vanishes at a point $\zeta_0$ near $bD$, because $T$ is constructed independently of $f_2$, it seems difficult to prove that $g_1$ is bounded excepted if we require that $g$ vanishes at $\zeta_0$ too. However, considering $g=f_2$, we easily see that in general this condition is not necessary when one wants to write $g$ as $g=g_1f_1+g_2f_2$ with $g_1$ and $g_2$ bounded for example. So the current in \cite{Maz2} probably does not give a good decomposition.

Actually, it appears that $f_2$ must be prioritized in the construction of the currents near a boundary point $\zeta_0$ such that $f_1(\zeta_0)=0$ and $f_2(z_0)\neq 0$ or more generally when $f_2$ is in some sense greater than $f_1$ and conversely. Following this idea, we construct two currents $T_1$ and $T_2$ such that $f_1T_1+f_2T_2=1$ on $D$. These currents are defined locally and using a suitable partition of unity we glue together the local currents and get a global current.

Let $\varepsilon _0$ be a small positive real number to be chosen later and let $\zeta_0$ be a point in $\overline D$. We distinguish three cases.
\par\medskip

If $\zeta_0$ belongs to $D_{-\varepsilon _0}$, we do not need to be careful.  
Using Weierstrass' preparation theorem when $\zeta_0$ belongs to $X_1$, we write $f_1=u_{0,1} P_{0,1}$ where $u_{0,1}$ is a non vanishing holomorphic function in a neighborhood ${\cal U}_0\subset D_{-\frac{\varepsilon _0}2}$ of $\zeta_0$ and $P_{0,1}(\zeta)=\zeta^{i_{0,1}}_2+ \zeta^{i_{0,1}-1}_2 a_{0,1}^{(1)}(\zeta_1) +\ldots + a_{0,1}^{(i_{0,1})}(\zeta_1)$, $a_{0,1}^{(k)}$ holomorphic on ${\cal U}_0$ for all $k$. If $\zeta_0$ does not belong to $X_1$, we set $P_{0,1}=1$, $i_{0,1}=0$, $u_{0,1}=f_1$ and we still have $f_1=u_{0,1} P_{0,1}$ with $u_{0,1}$ which does not vanish on some neighborhood ${\cal U}_0$ of $\zeta_0$. \\
For a smooth $(2,2)$-form $\varphi $ compactly supported in ${\cal U}_0$ we set
\begin{eqnarray*}
\langle T_{0,1}, \varphi\rangle&=&\frac1{c_0}\int_{{\cal U}_0} \frac{\overline{P_1(\zeta)}}{f_1(\zeta)}\diffp{^{i_{0,1}} \varphi }{\overline\zeta_2^{i_{0,1}}}(\zeta),\\
\langle T_{0,2}, \varphi\rangle&=&0,
\end{eqnarray*}
where $c_0$ is a suitable constant (see \cite{Maz2}). Integrating by parts we get $f_1T_{0,1}+f_2T_{0,2}=1$ on ${\cal U}_0$.

\par\medskip

If $\zeta_0$ belongs to $bD\setminus (X_1\cap X_2)$, without restriction we assume that $f_1(\zeta_0)\neq 0$. Let ${\cal U}_0$
 be a neighborhood of $\zeta_0$ such that $f_1$ does not vanish in ${\cal U}_0$. As in the previous case when $f_1(\zeta_0)\neq 0$, we set $P_{0,1}=1$, $i_{0,1}=0$, $u_{0,1}=f_1$ and for  any smooth $(2,2)$-form $\varphi $ compactly supported in $D\cap {\cal U}_0$ we put
 \begin{eqnarray*}
\langle T_{0,1}, \varphi\rangle&=&\frac1{c_0} \int_{{\cal U}_0} \frac{\overline{P_1(\zeta)}}{f_1(\zeta)}\diffp{^{i_{0,1}} \varphi }{\overline\zeta_2^{i_{0,1}}}(\zeta),\\
\langle T_{0,2}, \varphi\rangle&=&0.
\end{eqnarray*}
where as previously $c_0$ is a suitable constant. Again, we have $f_1T_{0,1}+f_2T_{0,2}=1$ on ${\cal U}_0\cap D$.

\par\medskip

If $\zeta_0$ belongs to $X_1\cap X_2\cap bD$, as in \cite{Al-Maz}, we cover a neighborhood ${\cal U}_0$ of $\zeta_0$ by a family of polydiscs $\pk1\zjk$, $j\in\nn$ and $k\in\{1,\ldots, n_j\}$ such that :
\begin{enumerate}[(i)]
 \item \label{seqi} For all $j\in\nn$, and all $k\in\{1,\ldots, n_j\}$, $\zjk$ belongs to $bD_{-(1-c\kappa )^j\varepsilon _0}$.
 \item \label{seqii} For all $j\in\nn$, all $k,l\in \{1,\ldots, n_j\}$, $k\neq l$, we have $\delta(z_{j,k},z_{j,l})\geq c\kappa (1-c\kappa )^j\varepsilon _0$.
 \item \label{seqiii} For all $j\in\nn$, all $z\in bD_{-(1-c\kappa )^j\varepsilon _0}$, there exists $k\in\{1,\ldots, n_j\}$ such that $\delta(z,\zjk)<c\kappa (1-c\kappa )^j\varepsilon _0$,
 \item \label{propmax1} $D\cap{\cal U}_0$ is included in $ \cup_{j=0}^{+\infty} \cup_{k=1}^{n_j} \pk1\zjk$,
 \item \label{propmax2} there exists $M\in\nn$ such that for $z\in D\setminus D_{-\varepsilon _0}$, ${\cal P}_{4\kappa |\rho (z)|}(z)$ intersect at most $M$ Koranyi balls ${\cal P}_{ 4\kappa |\rho (\zjk)|}\left(\zjk\right)$.
\end{enumerate}
Such a family of polydiscs will be called a $\kappa$-covering.

We define on each polydisc $\pk1\zjk$ two currents $T_{0,1}^{(j,k)}$ and $T_{0,2}^{(j,k)}$ such that $f_1T_{0,1}^{(j,k)}+f_2T^{(j,k)}_{0,2}=1$ as follows. We denote by $\Delta_\xi(\varepsilon )$ the disc of center $\xi$ and radius $\varepsilon $ and by  $(\zeta_{0,1}^*,\zeta_{0,2}^*)$ the coordinates of $\zeta_0$ in the \ko basis at $\zjk$. In \cite{Al-Maz} were proved the next two propositions~:
\begin{proposition}\label{propinter}
 If $\kappa>0$ is small enough and if $\pk1\zjk\cap X_l\neq\emptyset$ then $|\zeta^*_{0,1}|\geq 2\kappa|\rho(\zjk)|$.
\end{proposition}
We assume $\kappa$ so small that Proposition \ref{propinter} holds for both $X_1$ and $X_2$ with the same $\kappa$. When $|\zeta^*_{0,1}|\geq 2\kappa|\rho(\zjk)|$ then $X_l$ can be parametrized as follows (see \cite{Al-Maz})~:
\begin{proposition}
 If $|\zeta^*_{0,1}|\geq 2\kappa|\rho(\zjk)|$, for $l=1$ and $l=2$, there exists $p_l$ functions $\alpha^{(j,k)}_{l,1},\ldots, \alpha_{l,p_l}^{(j,k)}$ holomorphic on $\Delta_0(2\kappa|\rho(\zjk)|)$, there exists $r>0$, not depending from $j$ nor from $k$, and there exists $u_l^{(j,k)}$ holomorphic on the ball of center $\zeta_0$ and radius $r$ such that~:
\begin{enumerate}[(i)]
 \item $\diffp{\alpha_{l,i}^{(j,k)}}{\zeta^*_1}$ is bounded on $\Delta_0(2\kappa|\rho(\zjk)|)$ uniformly \wrt $j$ and $k$,
\item for all $\zeta\in\pk2\zjk$, $f_l(\zeta)=u_l^{(j,k)}(\zeta)\prod_{i=1}^{p_l}(\zeta_2^*-\alpha_{l,i}^{(j,k)}(\zeta^*_1))$.
\end{enumerate}
\end{proposition}
Now we define $T^{(j,k)}_{0,1}$ and $T_{0,2}^{(j,k)}$ with the following settings.\\
If $|z^*_{0,1}|<2\kappa|\rho(\zjk)|$ we set for $l=1$ and $l=2$~:
\begin{eqnarray*}
 I_l^{(j,k)}&:=&\emptyset;\\
 i_l^{(j,k)}&:=&0;\\
 P^{(j,k)}_l(\zeta)&:=&1.
\end{eqnarray*}
If $|z^*_{0,1}|\geq2\kappa|\rho(\zjk)|$ we set for $l=1$ and $l=2$~:
\begin{eqnarray*}
I_l^{(j,k)}&:=&\{i,\ \exists z^*_1\in\cc,\ |z^*_1|<\kappa |\rho (\zjk)|\text{ and } |\alpha_{l,i}^{(j,k)}(z_1^*)|<(2\kappa |\rho (\zjk)|)^{\frac12}\};\\
  i_l^{(j,k)}&:=&\#I_l^{(j,k)},\text{ the cardinal of } I_l^{(j,k)};\\
 P^{(j,k)}_l(\zeta)&:=&\prod_{i\in I_l^{(j,k)}}\left(\zeta^*_2-\alpha_{i,l}^{(j,k)}(\zeta_1^*)\right).
\end{eqnarray*}
In both case we set 
\begin{eqnarray*}
 {\cal U}_{1}^{(j,k)}&:=&\left\{\zeta\in \pk1\zjk,\ \left|\frac{f_1(\zeta)\rho(\zjk)^{i_1^{(j,k)}}}{P^{(j,k)}_1(\zeta)}\right|>\frac13\left|\frac{f_2(\zeta)\rho(\zjk)^{i_2^{(j,k)}}}{P^{(j,k)}_2(\zeta)}\right|\right\},\\
 {\cal U}_{2}^{(j,k)}&:=&\left\{\zeta\in \pk1\zjk,\frac23 \left|\frac{f_2(\zeta)\rho(\zjk)^{i_2^{(j,k)}}}{P^{(j,k)}_2(\zeta)}\right|>\left|\frac{f_1(\zeta)\rho(\zjk)^{i_1^{(j,k)}}}{P^{(j,k)}_1(\zeta)}\right|\right\},
\end{eqnarray*}
so that $\pk1\zjk={\cal U}_{1}^{(j,k)}\cup {\cal U}_{2}^{(j,k)}$.\\
For $l=1,2$ and for a smooth $(2,2)$-form $\varphi$ compactly supported in ${\cal U}_{l}^{(j,k)}$ we set
$$\langle T^{(j,k)}_{0,l}, \varphi \rangle :=\int_{\cc^2} \frac{\overline{P^{(j,k)}_l(\zeta)}}{f_l(\zeta)} \diffp{^{i_{l}^{(j,k)}}\varphi }{\overline{\zeta^*_2}^{i_l^{(j,k)}}}(\zeta).$$
Integrating $i_l^{(j,k)}$-times by parts, we get $f_lT_{0,l}^{(j,k)}=c_l^{(j,k)}$ on ${\cal U}_{l}^{(j,k)}$ 
where $c_l^{(j,k)}$ is an integer bounded by $i_{l}^{(j,k)}!$
(see \cite{Maz2}).
\par\medskip
Now we glue together the currents $T_{0,l}^{(j,k)}$ in order to define the current $T_{0,l}$, $l=1$, $2$, such that $f_1T_{0,1}+f_2T_{0,2}=1$ on $D\cap{\cal U}_0$. 
Let $(\tilde \chi_{j,k})_{\over{j\in\nn}{k\in\{1,\ldots,n_j\}}}$ be a partition of unity subordinated to the covering $(\pk1\zjk)_{\over{j\in\nn}{k\in\{1,\ldots,n_j\}}}$ of ${\cal U}_0$. Without restriction, we assume that $\left|\diffp{^{\alpha+\beta+\overline\alpha+\overline\beta} \tilde \chi _{j,k}}{{\zeta^*_1}^\alpha\partial {\zeta^*_2}^\beta\partial \overline{\zeta^*_1}^{\overline\alpha}\partial \overline{\zeta^*_2}^{\overline\beta}}(\zeta)\right|\leqs \frac1{|\rho (\zjk)|^{\alpha+\overline\alpha+\frac{\beta+\overline\beta}2}}.$ Let also $\chi $ be a smooth function on $\cc²\setminus\{0\}$  such that $\chi (z_1,z_2)=1$ if $|z_1|>\frac23 |z_2|$ and $\chi (z_1,z_2)=0$ if $|z_1|<\frac13|z_2|$ and let us define 
\begin{eqnarray*}
 \chi_{1}^{(j,k)}(\zeta)&=&\tilde \chi_{j,k}(\zeta)\cdot\chi\left(\frac{f_1(\zeta)\rho(\zjk)^{i_1^{(j,k)}}}{P^{(j,k)}_1(\zeta)},\frac{f_2(\zeta)\rho(\zjk)^{i_2^{(j,k)}}}{P^{(j,k)}_2(\zeta)} \right),\\
\chi_{2}^{(j,k)}(\zeta)&=&\tilde \chi_{j,k}(\zeta)\cdot\left(1-\chi\left(\frac{f_1(\zeta)\rho(\zjk)^{i_1^{(j,k)}}}{P^{(j,k)}_1(\zeta)},\frac{f_2(\zeta)\rho(\zjk)^{i_2^{(j,k)}}}{P^{(j,k)}_2(\zeta)} \right)\right).\\
\end{eqnarray*}
For $l=1$ and $l=2$, the support of $\chi_l^{(j,k)}$ is included in ${\cal U}^{(j,k)}_l$ so we can put
$$T_{0,l}=\sum_{\over{j\in\nn}{k\in\{1,\ldots, n_j\}}}\frac1{c^{(j,k)}_l} \chi_{l}^{(j,k)}T_{0,l}^{(j,k)}$$
and we have $f_1T_{0,1}+f_2T_{0,2}=1$ on ${\cal U}_0\cap D$.
\par\medskip

Now for all $\zeta_0\in bD \cup \overline{D_{-\varepsilon_0}}$ we have constructed a neighborhood ${\cal U}_0$ of $\zeta_0$ and two currents $T_{0,1}$ and $T_{0,2}$ such that $f_1T_{0,1}+f_2T_{0,2}=1$ on ${\cal U}_0\cap D$. If $\varepsilon_0>0$ is sufficiently small, we can cover $\overline D$ by finitely many open sets ${\cal U}_1,\ldots, {\cal U}_n$. Let $\chi_1,\ldots,\chi_n$ be a partition of unity subordinated to this family of open sets and $T_{1,1},\ldots, T_{k,1}$ and $T_{1,2},\ldots, T_{n,2}$ be the corresponding currents defined on ${\cal U}_1,\ldots, {\cal U}_n$. We glue together this current and we set
$$T_1=\sum_{j=1}^k \chi_j T_{j,1}\text{ and } T_2=\sum_{j=1}^n \chi_j T_{j,2},$$
so that $f_1T_1+f_2T_2=1$ on $D$. Moreover $T_1$ and $T_2$ are currents supported in $\overline D$ thus they have a finite order $k_2$ and we can apply $T_1$ and $T_2$ to function of class $C^{k_2}$ with support in $\overline D$. This gives $k_2$ from Theorem \ref{main_result}.

\section{The division formula}\label{section4}
In this part, given any two currents $T_1$ and $T_2$ of order $k_2$ such that $f_1T_1+f_2T_2=1$,  
assuming that $g$ is a holomorphic function on $D$ which belonging  to the ideal generated by $f_1$ and $f_2$, and which can be written as $g=\tilde g_1f_1+\tilde g_2f_2$, where $\tilde g_1$ and $\tilde g_2$ are two $C^\infty$-smooth functions on $D$ such that $|\rho|^N\tilde g_1$ and $|\rho |^N\tilde g_2$ vanish to order $k_2$ on $bD$ for some $N\in\nn$ sufficiently big,  we write $g$ as  $g= g_1f_1+ g_2f_2$ with $g_1$ and $g_2$ holomorphic on $D$. We point out that the formula we will get is valid for any $T_1$ and $T_2$ of order $k_2$ such that $f_1T_1+f_2T_2=1$.\\
Under our assumptions, for $k=1$ and $k=2$ and all fixed $z\in D$,  $\tilde g_1P^{N,k}(\cdot, z)$ and $\tilde g_2P^{N,k}(\cdot, z)$ can be extended by zero outside $D$ and are of class $C^{k_2}$ on $\cc$. So we can apply $T_1$ and $T_2$ to $\tilde g_1P^{N,k}(\cdot, z)$ and $\tilde g_2P^{N,k}(\cdot, z)$.
Now we construct a division formula.

For $l=1,2$, we denote by $b_l=b_{l,1}d\zeta_1+b_{l,2}d\zeta_2$ a $(1,0)$-form such that $f_l(z)-f_l(\zeta)=\sum_{i=1,2} b_{l,i}(\zeta,z)(z_i-\zeta_i)$. For the estimates, we will take $b_{l,i}(\zeta,z)=\int_0^1\diffp{f_l}{\zeta_i}(\zeta+t(z-\zeta))dt$, but this is not necessary to get a division formula.\\
From Theorem \ref{BA}, we have for all $z\in D$~:
$$g(z)=\int_Dg(\zeta)P^{N,2}(\zeta,z)$$
and since $g=\tilde g_1f_1+\tilde g_2f_2$
\begin{eqnarray}
 g(z)&=&\nonumber
 f_1(z)\int_D \tilde g_1(\zeta) P^{N,2}(\zeta,z)+f_2(z)\int_D \tilde g_2(\zeta) P^{N,2}(\zeta,z)\\
 & &+
 \int_D \tilde g_1(\zeta)\left(f_1(\zeta)-f_1(z)\right) P^{N,2}(\zeta,z)+\int_D \tilde g_2(\zeta)\left(f_2(\zeta)-f_2(z)\right) P^{N,2}(\zeta,z).\label{eq1}
\end{eqnarray}
Now from \cite{Maz1}, Lemma 3.4, there exists $\tilde c_{N,2}$ such that
\begin{eqnarray*}
\left(f_1(\zeta)-f_1(z)\right)P^{N,2}(\zeta,z)&=&\tilde c_{N,2} b_1(\zeta,z)\wedge\overline\partial P^{N,1}(\zeta,z)
\end{eqnarray*}
and since by assumption $\tilde g_1P^{N,1}$ vanishes on $bD$, Stokes' Theorem yields
\begin{eqnarray}
 \int_D \tilde g_1(\zeta)\left(f_1(\zeta)-f_1(z)\right) P^{N,2}(\zeta,z)&=&{\tilde c_{N,2}}\int_D\overline\partial \tilde g_1(\zeta)\wedge b_1(\zeta,z)\wedge P^{N,1}(\zeta,z).\label{eq2}
 \end{eqnarray}
We now use the fact that $f_1T_1+f_2T_2=1$ in order to rewrite the former integral :
\begin{eqnarray}
 \lefteqn{\int_D\overline\partial \tilde g_1(\zeta)\wedge b_1(\zeta,z)\wedge P^{N,1}(\zeta,z)}\nonumber\\
 &=&\langle f_1T_1+f_2T_2,\overline\partial \tilde g_1\wedge b_1(\cdot,z)\wedge P^{N,1}(\cdot,z)\rangle\nonumber\\
 &=&\langle f_1T_1,\overline\partial \tilde g_1\wedge b_1(\cdot,z)\wedge P^{N,1}(\cdot,z)\rangle+f_2(z)\langle T_2,\overline\partial \tilde g_1\wedge b_1(\cdot,z)\wedge P^{N,1}(\cdot,z)\rangle\nonumber\\
 &&+\langle T_2,\left(f_2-f_2(z)\right)\overline\partial \tilde g_1\wedge b_1(\cdot,z)\wedge P^{N,1}(\cdot,z)\rangle.\label{eq3}
\end{eqnarray}
Again from \cite{Maz1}, Lemma 3.4, there exists $\tilde c_{N,1}$ such that
\begin{eqnarray*}
\lefteqn{\left(f_2(\zeta)-f_2(z)\right) b_1(\zeta,z)\wedge P^{N,1}(\zeta,z) -
\left(f_1(\zeta)-f_1(z)\right) b_2(\zeta,z)\wedge P^{N,1}(\zeta,z)}\\
&&\hskip 200pt= \tilde c_{N,1} b_1(\zeta,z)\wedge b_2(\zeta,z)\wedge \overline\partial P^{N,0}(\zeta,z).
\end{eqnarray*}
So 
\begin{eqnarray}
 \lefteqn{\langle T_2,\left(f_2-f_2(z)\right)\overline\partial  \tilde g_1\wedge b_1(\cdot,z)\wedge P^{N,1}(\cdot,z)\rangle}\nonumber\\
 &=& -f_1(z)\langle T_2, \overline\partial \tilde g_1\wedge b_2(\cdot,z)\wedge P^{N,1}(\cdot,z)\rangle+
 \langle T_2, f_1\overline\partial  \tilde g_1\wedge b_2(\cdot,z)\wedge P^{N,1}(\cdot,z)\rangle\nonumber\\
 &&+
 \tilde c_{N,1} 
\langle T_2,\overline\partial  \tilde g_1\wedge b_1(\cdot,z)\wedge b_2(\cdot,z)\wedge \overline\partial P^{N,0}(\cdot,z)\rangle\label{eq4}
\end{eqnarray}
We plug  together (\ref{eq2}), (\ref{eq3}) and (\ref{eq4}) and their analogue for $\int_D g_2(\zeta)\left(f_2(\zeta)-f_2(z)\right) P^{N,2}(\zeta,z)$ in  (\ref{eq1}) and we get
\begin{eqnarray}
 g(z)&=&\nonumber
 f_1(z)\int_D \tilde  g_1(\zeta) P^{N,2}(\zeta,z)
 -\tilde c_{N,2} f_1(z)\langle T_2, \overline\partial  \tilde g_1\wedge b_2(\cdot,z)\wedge P^{N,1}(\cdot,z)\rangle\\
 &&
 +\tilde c_{N,2} f_2(z)\langle T_2,\overline\partial  \tilde g_1\wedge b_1(\cdot,z)\wedge P^{N,1}(\cdot,z)\rangle\nonumber\\
 &&+f_2(z)\int_D  \tilde g_2(\zeta) P^{N,2}(\zeta,z)-\tilde c_{N,2} f_2(z)\langle T_1, \overline\partial \tilde  g_2\wedge b_1(\cdot,z)\wedge P^{N,1}(\cdot,z)\rangle\nonumber\\
 &&+\tilde c_{N,2} f_1(z)\langle T_1,\overline\partial  \tilde g_2\wedge b_2(\cdot,z)\wedge P^{N,1}(\cdot,z)\rangle\nonumber\\
 &&+\tilde c_{N,2} \langle T_1,f_1\overline\partial  \tilde g_1\wedge b_1(\cdot,z)\wedge P^{N,1}(\cdot,z)\rangle+
 \tilde c_{N,2} \langle T_2, f_1\overline\partial  \tilde g_1\wedge b_2(\cdot,z)\wedge P^{N,1}(\cdot,z)\rangle\label{eq5}\\
 &&+\tilde c_{N,2} \langle T_2,f_2\overline\partial  \tilde g_2\wedge b_2(\cdot,z)\wedge P^{N,1}(\cdot,z)\rangle+
\tilde c_{N,2}  \langle T_1, f_2\overline\partial  \tilde g_2\wedge b_1(\cdot,z)\wedge P^{N,1}(\cdot,z)\rangle\label{eq6}\\
 &&+ \tilde c_{N,2} \tilde c_{N,1} 
\langle \overline\partial  \tilde g_1\wedge T_2-\overline\partial  \tilde g_2\wedge T_1 ,b_1(\cdot,z)\wedge b_2(\cdot,z)\wedge \overline\partial P^{N,0}(\cdot,z)\rangle\nonumber
 \end{eqnarray}
 Now since $\overline\partial g=f_1\overline\partial \tilde  g_1+f_2\overline\partial  \tilde g_2=0$, the line (\ref{eq5}) and (\ref{eq6}) vanish. Therefore in order to get our division formula, it suffices to prove that  $\overline\partial(\overline \partial  \tilde g_1\wedge T_2-\overline\partial  \tilde g_2\wedge T_1)=0$.\\
When $X_1\cap X_2$ is not a complete intersection and when assumption $(iv)$ in the introduction is satisfied by $ \tilde g_1$ and $ \tilde g_2$, one can prove that $\overline \partial  \tilde g_1\wedge \overline\partial T_2=0$ and $\overline \partial  \tilde g_2\wedge \overline \partial T_1=0$.\\
When $X_1\cap X_2$ is a complete intersection, we prove that for any $\zeta_0\in D$ there exists a neighborhood ${\cal U}_0$ of $\zeta_0$ such that for all  $(2,1)$-form $\varphi$, smooth and supported in ${\cal U}_0$, we have 
$\langle\overline \partial \tilde g_1\wedge T_2-\overline\partial \tilde  g_2\wedge T_1,\overline\partial\varphi\rangle=0$.\\
Let $\zeta_0$ be a point in $D$. By assumption on $g$, there exists a neighborhood ${\cal U_0}$ of $\zeta_0$ and two holomorphic functions $\gamma_1$ and $\gamma_2$ such that $g=\gamma_1f_1+\gamma_2f_2$ on ${\cal U}_0$. 
We now use the following lemma from which we postpone the proof to the end of this section~:
\begin{lemma}\label{to_be_written}
 Let $f_1$ and $f_2$ be two holomorphic functions defined in a neighborhood of $0$ in $\cc^2$, $X_1=\{z,\ f_1(z)=0\}$ and $X_2=\{z,\ f_2(z)=0\}$. We assume that $X_1\cap X_2$ is a complete intersection and that $0$ belongs to $X_1\cap X_2$. Let $\varphi_1$ and $\varphi_2$ be two $C^\infty$-smooth functions such that $f_1\varphi_1=f_2\varphi_2$.\\
Then, $\frac{\varphi_1}{f_2}$ and $\frac{\varphi_2}{f_1}$ are $C^\infty$-smooth in a neighborhood of $0$.
\end{lemma}
Lemma \ref{to_be_written} implies that  the function $\psi=\frac{ \tilde g_1-\gamma_1}{f_2}=\frac{\gamma_2- \tilde g_2}{f_1}$ is smooth on a perhaps smaller neighborhood of $\zeta_0$ still denoted by ${\cal U}_0$. Thus
\begin{eqnarray*}
\langle\overline \partial  \tilde g_1\wedge T_2-\overline\partial \tilde  g_2\wedge T_1,\overline\partial\varphi\rangle&=&
\langle\overline \partial ( \tilde g_1-\gamma_1)\wedge T_2+\overline\partial (\gamma_2- \tilde g_2)\wedge T_1,\overline\partial\varphi\rangle\\
&=&\langle\overline \partial (f_2\psi)\wedge T_2+\overline\partial (f_1\psi )\wedge T_1,\overline\partial\varphi\rangle\\
&=&\langle f_2 T_2+f_1 T_1,\overline\partial\psi\wedge \overline\partial\varphi\rangle\\
&=&\int_{{\cal U}_0}\overline\partial\psi\wedge \overline\partial\varphi
\end{eqnarray*}
and since $\varphi$ is supported in ${\cal U}_0$ we have $\int_{{\cal U}_0}\overline\partial\psi\wedge \overline\partial\varphi=-\int_{{\cal U}_0} d(\varphi\overline\partial \psi)=0$ and so 
$$\langle\overline \partial  \tilde g_1\wedge T_2-\overline\partial  \tilde g_2\wedge T_1,\overline\partial\varphi\rangle=0.$$
Now we set
\begin{eqnarray*}
\lefteqn{g_1(z)=\int_D  \tilde g_1(\zeta) P^{N,2}(\zeta,z)}\\
&& +\tilde c_{N,2}\left(\langle T_1,\overline\partial  \tilde g_2\wedge b_2(\cdot,z)\wedge P^{N,1}(\cdot,z)\rangle-
 \langle T_2, \overline\partial  \tilde g_1\wedge b_2(\cdot,z)\wedge P^{N,1}(\cdot,z)\rangle\right)\\
\lefteqn{g_2(z)=\int_D  \tilde g_2(\zeta) P^{N,2}(\zeta,z)}\\
&& +\tilde c_{N,2}\left(\langle T_2,\overline\partial  \tilde g_1\wedge b_1(\cdot,z)\wedge P^{N,1}(\cdot,z)\rangle-
 \langle T_1, \overline\partial  \tilde g_2\wedge b_1(\cdot,z)\wedge P^{N,1}(\cdot,z)\rangle\right)
\end{eqnarray*}
and we have
$$g=g_1f_1+g_2f_2$$
with $g_1$ and $g_2$ holomorphic on $D$. We notice that if $ \tilde g_1$ and $ \tilde g_2$ are already holomorphic functions then $g_1=\tilde g_1$ and $g_2=\tilde g_2$.\\[10pt]
{\it Proof of Lemma \ref{to_be_written}~:} 
Maybe after a unitary change of coordinates, we can assume that for $l=1,2$, the function $f_l$ is given by $f_l(z,w)= z^{k_l}+a^{(l)}_{1}(w)z^{k_l-1}+\ldots+a_{k_l}^{(l)}(w)$ where $a^{(l)}_{1},\ldots, a^{(l)}_{k_l}$ are holomorphic near $0$ and vanish at $0$. Moreover, since the intersection $X_1\cap X_2$ is transverse, $P_1$ and $P_2$ are relatively prime. Thus there exists two polynomials $\alpha_1$ and $\alpha_2$ with holomorphic coefficients in $w$ and a function $\beta$ of $w$ not identically zero such that 
$$\alpha_1(z,w)f_1(z,w)+\alpha_2(z,w)f_2(z,w)=\beta(w).$$
Multiplying this equality by $\varphi_1$ we get
$$ f_2(\alpha_1\varphi_2+\alpha_2\varphi_1)=\beta\varphi_1.$$
We now prove that $\beta$ divides the function $\psi:=\alpha_1\varphi_2+\alpha_2\varphi_1$.\\
Since $\beta$ is not identically zero, there exists $k\in\nn$ such that $\beta(w)=w^k \gamma(w)$ where $\gamma(0)\neq 0$.\\
For all $j\in\nn$ we have 
\begin{eqnarray}
f_2(z,w)\diffp{^j\psi}{\overline w^j}(z,w)&=&\beta(w)\diffp{\varphi_1}{\overline w^j}(z,w) \label{eq7}
\end{eqnarray}
and for $w=0$ and all $z$ we thus get $\diffp{^j\psi}{\overline w^j}(z,0)=0$.\\
By induction we then deduce from (\ref{eq7}) that $\diffp{^{i+j}\psi}{w^i\partial \overline w^j}(z,0)=0$ for all $i\in\{0,\ldots, k-1\}$ and all $j\in\nn$.
For any integer $n\geq k$ we therefore can write for all $z$ and all $w$
$$\frac{\psi(z,w)}{w^k}=\sum_{\over{k\leq i+j \leq n}{i\geq k}} w^{i-k} \overline w^j \diffp{^{i+j} \psi}{w^i\partial \overline w^j}(z,0)
+\sum_{i+j=n+1} w^{i-k}\overline w^j\int_0^1\diffp{^{n+1}\psi}{w^i\partial \overline w^j}(z,tw)dt.$$
Now, it is easy to check by induction that the function $w\mapsto \frac{\overline w^{i+j}}{w^i}$ is of class $C^{j-1}$ for all positive integer $j$ and all non negative integer $i$. This implies that $\frac{\psi(z,w)}{w^k}$ is of class $C^{n}$ for all positive integer $n$ and therefore $\frac{\varphi_1}{f_2}=\frac{\psi}{\beta}$ is of class $C^\infty$.\qed

\section{Proof of the main result}\label{section5}
In order to prove Theorem \ref{main_result}, for any $k$ and $l$ in $\{1,2\}$ and any $q\in[1,+\infty]$, we have to prove that if $h$ is a smooth function such that, for all non-negative integers $\alpha$ and $\beta$, $\left|\diffp{^{\alpha+\beta} \tilde h}{\overline{\eta_\zeta}^\alpha\partial\overline{v_\zeta}^\beta}\right| |\rho|^{\alpha+\frac\beta2}$ belongs to $L^q(D)$,  then the function
$$z\mapsto \langle T_l,\overline\partial h\wedge b_{k}(\cdot,z)\wedge P^{N,1}(\cdot,z)\rangle$$
belongs to $L^q(D)$ if $q<\infty$ and to $BMO(D)$ if $q=+\infty$.

As usually, the main difficulty occurs when $z$ is near $bD$ and when we integrate for $\zeta$ near $z$. Moreover, the only interesting case here is when, in addition, $z$ is near a point $\zeta_0\in bD\cap X_1\cap X_2$ and we only consider that case.\\
 We use the same notation as in section \ref{construction} and assume that $z$ belongs to the neighborhood ${\cal U}_0$ of a point $\zeta_0\in bD\cap X_1\cap X_2$ which was used during the construction of the currents. Moreover, we assume that the $\ko$ basis at $\zeta_0$ is the canonical basis of $\cc^2$ and that $\zeta_0$ is the origin of $\cc^2$. We will need upper bound of $\frac{P_l^{(j,k)}}{f_l}\diffp{^{\alpha+\beta}f_l} {{\zeta^*_1}^\alpha \partial {\zeta^*_2}^ \beta}$ in order to estimate $\frac{P_l^{(j,k)}}{f_l}b_m$ and the derivatives of $\chi_l^{(j,k)}$. We begin with the following lemma~:
\begin{lemma}\label{lem1}
 For all $j\in\nn$, all $k\in\{1,\ldots, n_j\},$ all $\alpha$ and $\beta$ in $\nn$, $l=1,2$, all $\zeta$ in $\pk1\zjk$ and all
 $\tilde \zeta\in\cc^2$ such that $|\tilde \zeta^*_1|<2\kappa|\rho(\zjk)|$ and $|\tilde \zeta^*_2|<(4\kappa|\rho(\zjk)|)^{\frac12}$, 
 we have uniformly with respect to $j,k,l,\zeta$ and $\tilde\zeta$
$$\left|
\frac{P^{(j,k)}_l(\zeta)}{f_l(\zeta)}\diffp{^{\alpha+\beta}}{{\zeta^*_1}^\alpha\partial {\zeta_2^*}^\beta} \left(\frac{f_l(\tilde \zeta)}{P_l^{(j,k)}(\tilde \zeta)}\right)
\right|\leqs |\rho(\zjk)|^{-\alpha-\frac\beta2}.$$
\end{lemma}
\pr We denote by $(\zeta_{0,1}^*,\zeta_{0,2}^*)$ the coordinates of $\zeta_0$ in the \ko coordinates at $z_{j,k}$. The definition of $P^{(j,k)}_l$ forces us to distinguish three cases~:
\par\medskip
{\it First case~:} If $|\zeta^*_{0,1}|<2\kappa|\rho(\zjk)|$ and $|\zeta^*_{0,2}|<(6\kappa|\rho(\zjk)|)^{\frac12}$, then $\delta(\zjk,\zeta_0)\leq 6\kappa|\rho(\zjk)|$ and thus for all $\tilde \zeta\in\pk6\zjk$, $\delta(\tilde \zeta,\zeta_0)\leqs |\rho(\zjk)|$.\\
For all $\varepsilon>0$ and all $\tilde \zeta\in{\cal P}_\varepsilon(\zeta_0)$, it is easy to see that $|f_l(\tilde \zeta)|\leqs \varepsilon^{\frac{p_l}2}$. Therefore, Cauchy's inequalities give
$$\left|\diffp{^{\alpha+\beta}f_l}{{\zeta^*_1}^\alpha\partial{\zeta^*_2}^\beta}(\tilde \zeta)\right|\leqs |\rho(\zjk)|^{\frac{p_l}2-\alpha-\frac\beta2}$$ for all $\zeta \in\pk4\zjk$. Moreover, since $|\zeta_{0,1}^*|<2\kappa|\rho(\zjk)|$, on the one hand $P^{(j,k)}_l=1$, and on the other hand $\pk1\zjk\cap X_l=\emptyset$ (see Proposition \ref{propinter}) which implies that $|f_l(\zeta)|\geqs |\rho(\zjk)|^{\frac{p_l}2}$ for all $\zeta\in\pk1\zjk$. Therefore
$\left|
\frac{P^{(j,k)}_l(\zeta)}{f_l(\zeta)}\diffp{^{\alpha+\beta}}{{\zeta^*_1}^\alpha\partial {\zeta_2^*}^\beta} \left(\frac{f_l(\tilde\zeta)}{P_l^{(j,k)}(\tilde \zeta)}\right)
\right|\leqs |\rho(\zjk)|^{-\alpha-\frac\beta2}.$
\par\medskip
{\it Second case~:} If $|\zeta^*_{0,1}|<2\kappa|\rho(\zjk)|$ and $|\zeta_{0,2}^*|\geq (6\kappa|\rho(\zjk)|)^{\frac12}$, we set $a(\zjk)=\diffp{\rho}{\zeta_1}(\zjk)$, $b(\zjk)=\diffp{\rho}{\zeta_2}(\zjk)$ and 
$$P(\zjk)=\frac1{\sqrt{|a(\zjk)|^2+|b(\zjk)|^2}}\left(\begin{array}{cc} a(\zjk)&b(\zjk)\\ -\overline{b(\zjk)}&\overline{a(\zjk)}\end{array}\right).$$
Then we have $\zeta^*=P(\zjk)(\zeta-\zjk)$. Moreover $b(\zjk)$ tends to 0 when $\zjk$ goes to $\zeta_0$, that is if ${\cal U}_0$ is sufficiently small.\\
For $\tilde \zeta\in\pk5\zjk$, if ${\cal U}_0$ is sufficiently small
\begin{eqnarray*}
 |\tilde \zeta_2|&\geq& \frac1{\sqrt{|a(\zjk)|^2+|b(\zjk)|^2}} (|a(\zjk)||\zeta_{0,2}^*|-|b(\zjk)||\zeta^*_{0,1}|-|b(\zjk)||\tilde \zeta^*_1|-|a(\zjk)||\tilde \zeta^*_2|)\\
&\geqs& |\zeta^*_{0,2}|.
\end{eqnarray*}
We also trivially have $|\tilde\zeta_2|\leqs |\zeta^*_{0,2}|
$ and so $|\tilde\zeta_2|\eqs |\zeta^*_{0,2}|$. Analogously we have $|\zeta_2|\eqs|\zeta^*_{0,2}|$ for $\zeta\in\pk1\zjk$.\\
On the other hand 
\begin{eqnarray*}
|\tilde \zeta_1|&\leq& \frac1{\sqrt{|a(\zjk)|^2+|b(\zjk)|^2}} \left(|a(\zjk)|( |\zeta_{0,1}^*|+|\tilde\zeta^*_1|) + |b(\zjk)|(|\zeta^*_{0,2}|+|\tilde\zeta^*_2|) \right)\\
&\leq &2\kappa|\rho(\zjk)|+|b(\zjk)|(|\zeta^*_{0,2}|+|\rho(\zjk)|^{\frac12})\\
&\leq &c |\zeta^*_{0,2}| 
\end{eqnarray*}
where $c$ does not depend from $\zjk$ nor from $\tilde \zeta$ and is arbitrarily small provided ${\cal U}_0$ is small enough. We also have $|\zeta_1|\leq c|\zeta^*_{0,2}|$ for $\zeta\in\pk1\zjk$.\\
Now let $\alpha\in \cc$ be such that $f_{l}(\zeta_1,\alpha)=0$. Since the intersection $X_l\cap bD$ is transverse, there exists a positive constant $C$ not depending from $\tilde \zeta$, $\alpha$, $j$ nor $k$ such that $|\alpha|\leq C|\tilde \zeta_1|$.\\
Therefore if $c$ is small enough, $|\alpha|\leq \frac12|\tilde \zeta_2|$. This yields
\begin{eqnarray*}
 |f_l(\tilde \zeta)|&\eqs& \prod_{\alpha / f_l(\tilde \zeta_1,\alpha)=0} |\tilde \zeta_2-\alpha|\\
&\eqs& 
|\zeta^*_{0,2}|
^{p_l}.
\end{eqnarray*}
Analogously we have $|f_l(\zeta)|\eqs|\zeta^*_{0,2}|^{p_l}$ for $\zeta\in\pk1\zjk$. 
Cauchy's inequalities then give for all $\tilde\zeta\in\pk4\zjk$ $$\left|
\diffp{^{\alpha+\beta}f_l}{{\zeta^*_1}^\alpha\partial {\zeta_2^*}^\beta}(\zeta)
\right|\leqs |\zeta^*_{0,2}|^{p_l} |\rho(\zjk)|^{-\alpha-\frac\beta2},$$
and since $P^{(j,k)}_l=1$ when $|\zeta^*_{0,1}|\leq 2\kappa|\rho(\zjk)|$, we are done in this case.
\par\medskip
{\it Third case~:} If $|\zeta^*_{0,1}|>2\kappa|\rho(\zjk)|$, there exists a family of parametrization $\alpha^{(j,k)}_{l,i}$, $i=1,\ldots, p_l,$ given by Proposition \ref{propinter} such that $\left|\diffp{^n\alpha^{(j,k)}_{l,i}}{{\zeta^*_1}^n}(\zeta^*_1)\right|\leqs |\rho(\zjk)|^{1-n}$ for all $\zeta_1^*\in\Delta_0(2\kappa|\rho(\zjk)|)$. Moreover is this case, we actually seek an upper bound for
$$\frac1{\prod_{i\notin I_l^{(j,k)}} \left(\zeta^*_2-\alpha^{(j,k)}_{l,i}(\zeta^*_1)\right)}\diffp{^{\alpha+\beta} }{{\zeta^*_1}^\alpha\partial {\zeta^*_2}^\beta}\left({\prod_{i\notin I_l^{(j,k)}} \left(\tilde \zeta^*_2-\alpha^{(j,k)}_{l,i}(\tilde \zeta^*_1)\right)}\right) .$$
We fix $i$ in $\{1,\ldots, p_l\}\setminus I_l^{(j,k)}$ and $\tilde\zeta$ such that $|\tilde \zeta^*_1|<2\kappa |\rho (\zjk)|$ and $|\tilde \zeta^*_2|<(4\kappa |\rho (\zjk)|)^{\frac12}$.
  If $|\alpha_{l,i}^{(j,k)}(\tilde\zeta_1^*)|\leq (6\kappa |\rho (\zjk)|)^{\frac12}$, we have $|\tilde \zeta^*_2-\alpha_{l,i}^{(j,k)}(\tilde \zeta^*_1)|\leqs (\kappa |\rho (\zjk)|)^{\frac12}$. On the other hand, by definition of $I^{(j,k)}_l$, for all $\zeta^*_1\in\Delta_0(2\kappa |\rho (\zjk)|)$, we have $|\alpha_{l,i}^{(j,k)}(\zeta_1^*)|\geq (2\kappa |\rho(z_{j,k})|)^{\frac12}$. Therefore, for all $\zeta\in\pk1\zjk$ we have
 \begin{eqnarray}
 \frac{ |\tilde \zeta^*_2-\alpha_{l,i}^{(j,k)}(\tilde \zeta_1^*)|}{|\zeta^*_2-\alpha_{l,i}^{(j,k)}(\zeta^*_1)|}&\leqs&1\label{eq50}
 \end{eqnarray}
uniformly with respect to $\zeta,\tilde\zeta$ and $\zjk$.\\
 If now $|\alpha_{l,i}^{(j,k)}(\tilde\zeta_1^*)|\geq (6\kappa |\rho (\zjk)|)^{\frac12}$, we have $|\tilde \zeta^*_2-\alpha_{l,i}^{(j,k)}(\tilde \zeta^*_1)|\leqs |\alpha_{l,i}^{(j,k)}(\tilde\zeta_1^*)|$.\\
 If ${\cal U}_0$ is sufficiently small, Proposition \ref{propII.0.2} then yields $|\alpha_{l,i}^{(j,k)}(\zeta^*_1)|\geqs (3\kappa |\rho (\zjk)|)^{\frac12} $ for all $\zeta\in\pk1\zjk$ and
 \begin{eqnarray}
 \frac{ |\tilde \zeta^*_2-\alpha_{l,i}^{(j,k)}(\tilde\zeta_1^*)|}{|\zeta^*_2-\alpha_{l,i}^{(j,k)}(\zeta^*_1)|}&\leqs&1\label{eq51}
 \end{eqnarray}
uniformly with respect to $\zeta,\tilde\zeta$ and $\zjk$.\\
From proposition \ref{propII.0.2} we also have  
\begin{eqnarray}
 \frac{ 1}{|\zeta^*_2-\alpha_{l,i}^{(j,k)}(\zeta^*_1)|}&\leqs&|\rho(\zjk)|^{-\frac12}\label{eq52}
 \end{eqnarray} for all $\zeta\in\pk1\zjk$ 
and $\left|\diffp{^\alpha\alpha_{l,i}^{(j,k)}}{\zeta_1^*}(\tilde\zeta^*_1)\right|\leqs |\rho(\zjk)|^{1-\alpha}$ for all $\tilde\zeta^*_1\in\Delta_0(2\kappa|\rho(\zjk)|)$ so
\begin{eqnarray}
 \frac{ 1}{|\zeta^*_2-\alpha_{l,i}^{(j,k)}(\zeta^*_1)|}\left|\diffp{^\alpha\alpha_{l,i}^{(j,k)}}{\zeta_1^*}(\tilde\zeta^*_1)\right|&\leqs&|\rho(\zjk)|^{-\alpha}.\label{eq53}
 \end{eqnarray}
Now the inequalities (\ref{eq50}), (\ref{eq51}), (\ref{eq52}) and (\ref{eq53}) yield the lemma.\qed
\par\medskip
Lemma \ref{lem1} gives us an  upper bound for the derivatives of $\chi_l^{(j,k)}$~:
\begin{corollary}\label{cor2}
 For all $j\in\nn$, all $k\in\{1,\ldots, n_j\},$ all $\alpha$ and $\beta$ in $\nn$, $l=1,2$ and all $\zeta\in\pk1\zjk$, we have uniformly with respect to $j,k,l$ and $\zeta$
$$\left|
\diffp{^{\alpha+\beta}\chi_l^{(j,k)}}{{\overline\zeta^*_1}^\alpha\partial {\overline\zeta_2^*}^\beta}(\zeta)
\right|\leqs |\rho(\zjk)|^{-\alpha-\frac\beta2}.$$
\end{corollary}
\pr Since by construction $\left|
\diffp{^{\alpha+\beta}\tilde \chi_{j,k}}{{\overline\zeta^*_1}^\alpha\partial {\overline\zeta_2^*}^\beta}(\zeta)
\right|\leqs |\rho(\zjk)|^{-\alpha-\frac\beta2}$, we only have to consider 
$\diffp{^{\alpha+\beta}}{{\zeta^*_1}^\alpha\partial {\zeta^*_2}^\beta}\chi\left(\frac{f_1(\zeta)}{P^{(j,k)}_1(\zeta)} \rzjk^{i_1^{(j,k)}},\frac{f_2(\zeta)}{P^{(j,k)}_2(\zeta)} \rzjk^{i_2^{(j,k)}} \right) $.\\
The derivative $\diffp{^{\gamma+\delta}\chi}{z_1^\gamma\partial z_2^\delta}(z_1,z_2)$ is bounded up to a uniform multiplicative constant by $\frac1{|z_1|^\gamma|z_2|^\delta}$ when $\frac13|z_2|<|z_1|<\frac23|z_2|$ and is zero otherwise.\\
Therefore, we can estimate $\left|
\diffp{^{\alpha+\beta}\chi_l^{(j,k)}}{\overline{\zeta^*_1}^\alpha\partial \overline{\zeta_2^*}^\beta}
\right|$ by a sum of products of $\left|\frac{P_l^{(j,k)}}{f_l}\diffp{^{\gamma+\delta}}{{\overline\zeta_1^*}^\gamma\partial{\overline\zeta^*_2}^\delta}\left(\frac{f_l}{P_l^{(j,k)}}\right)\right|$ where the sum of the $\gamma$'s equals $\alpha$ and the sum of the $\delta$'s equals $\beta$. Lemma \ref{lem1} then gives the wanted estimates.\qed
\begin{corollary}\label{cor4}
For any smooth function $h$, we can write 
$$\diffp{^{i_l^{(j,k)}}}{\overline{\zeta^*_2}^{i_l^{(j,k)}}} \left(\chi_l^{(j,k)}(\zeta) \overline\partial h(\zeta)\wedge P^{N,1}(\zeta,z)\right)=\psi_1^{(j,k,l)}(\zeta,z)d\zeta^*_1+\psi_2^{(j,k,l)}(\zeta,z)d\zeta^*_2$$
 with $\psi_1^{(j,k,l)}$ and $\psi_2^{(j,k,l)}$ two (0,2)-forms supported in  ${\cal U}^{(j,k)}_l$ satisfying for $\nabla_z$ a differential operator of order $1$ acting on $z$, uniformly with respect to $j,k,z$ and $\zeta\in{\cal U}^{(j,k)}_l$~:
\begin{eqnarray*}
 {\left|\psi_1^{(j,k,l)}(\zeta,z)\right|}
&\leqs& \rzjk^{-i_l^{(j,k)}-\frac52} \left(\frac{\rzjk}{\rzjk+|\rho(z)|+\delta(\zjk,z)}\right)^{N} \tilde h(\zeta),\\
{\left|\psi_2^{(j,k,l)}(\zeta,z)\right|}
&\leqs& \rzjk^{-i_l^{(j,k)}-2} \left(\frac{\rzjk}{\rzjk+|\rho(z)|+\delta(\zjk,z)}\right)^{N}\tilde h(\zeta),\\
{\left|\nabla_z{\psi_1^{(j,k,l)}}(\zeta,z)\right|}
&\leqs&  \rzjk^{-i_l^{(j,k)}-\frac72} \left(\frac{\rzjk}{\rzjk+|\rho(z)|+\delta(\zjk,z)}\right)^{N} \tilde h(\zeta),\\
{\left|\nabla_z{\psi_2^{(j,k,l)}}(\zeta,z)\right|}
&\leqs& \rzjk^{-i_l^{(j,k)}-3} \left(\frac{\rzjk}{\rzjk+|\rho(z)|+\delta(\zjk,z)}\right)^{N}\tilde h(\zeta),
\end{eqnarray*}
where $\tilde h(\zeta)=\max_{n\in\{0,\ldots, i_l^{(j,k)}\}} \left(\left|\diffp{^{n+1} h}{\overline{\zeta^*_2}^{n+1}}(\zeta)|\rho(\zeta)|^{\frac {n+1} 2} \right|, \left|\diffp{^{n+1} h}{\overline{\zeta^*_1}\partial \overline{\zeta^*_2}^{n}}(\zeta)|\rho(\zeta)|^{\frac {n} 2+1} \right|\right)$.
\end{corollary}
\pr Proposition \ref{estiBA} implies that $\diffp{^n}{\overline{\zeta^*_2}^n}P^{N,1}(\zeta,z)=\sum_{p,q=1,2} \tilde\psi^{(n,N)}_{p,q}(\zeta,z)d\zeta^*_p\wedge d\overline{\zeta_q^*}$ where 
$$|\tilde\psi^{n,N}_{p,q}(\zeta,z)|\leqs \left(\frac{|\rho(\zeta)|}{|\rho(\zeta)|+|\rho(z)|+\delta(\zeta,z)}\right)^N|\rho(\zeta)|^{-\frac1p-\frac1q-\frac n2}.$$ 
From proposition \ref{propII.0.2}, if $\kappa$ is small enough, we have for all $\zeta\in \pk1\zjk$, $\frac12\rzjk\leq|\rho(\zeta)|$ and thus, provided $\kappa$ is small enough~:
\begin{eqnarray*}
|\rho(\zeta)|+\delta(\zeta,z)
&\geq& \frac12\rzjk+\frac1{c_1} \delta(z,\zjk)-\delta(\zjk,\zeta)\\
&\geqs& \rzjk+\delta(z,\zjk)
\end{eqnarray*}
and so $|\tilde\psi^{n,N}_{p,q}(\zeta,z)|\leqs \left(\frac{|\rho(\zjk)|}{|\rho(\zjk)|+|\rho(z)|+\delta(\zjk,z)}\right)^N|\rho(\zjk)|^{-\frac1p-\frac1q-\frac n2}.$ This inequality and Corollary \ref{cor2} now yield the two first estimates. The two others can be shown in the same way.\qed\par\medskip
In order to estimate $\frac{\overline{P^{(j,k)}_l}}{f_l}b_k$, we need the following lemma~:
\begin{lemma}\label{lem2}
For all $j\in\nn$, all $k\in\{1,\ldots, n_j\},$ all $\alpha$ and $\beta$ in $\nn$, $l=1,2$ and all $\zeta\in\pk1\zjk$ we have uniformly with respect to $j,k,l$ and $\zeta$
$$\left|
\diffp{^{\alpha+\beta}}{{\zeta^*_1}^\alpha\partial {\zeta_2^*}^\beta} \left(
\prod_{i\in I_l^{(j,k)}} (\zeta^*_2-\alpha^{(j,k)}_{l,i}(\zeta^*_1))\right)
\right|\leqs |\rho(\zjk)|^{i^{(j,k)}_l-\alpha-\frac\beta2}.$$
\end{lemma}
\pr For $i\in I^{(j,k)}_{l}$, there exists $z^*_1\in \Delta_0(\kappa|\rho(\zjk)|)$ such that $|\alpha^{(j,k)}_{l,i}(z^*_1)|<2\kappa\rzjk^{\frac12}$. Since $\left|\diffp{\alpha_{l,i}^{(j,k)}}{\zeta_1^*}(\zeta_1^*)\right|$ is uniformly bounded on $\Delta_0(2\kappa\rzjk)$, for all $\zeta\in{\cal P}_{2\kappa\rzjk}(z_{j,k})$, we have $\prod_{i\in I^{(j,k)}_l}\left |\zeta^*_2-\alpha^{(j,k)}_{l,i}(\zeta^*_1)\right|\leqs \rzjk^{\frac{i_l^{(j,k)}}2}$. Cauchy's inequalities then give the results.\qed
\par\medskip
As a direct corollary of Lemma \ref{lem1} and \ref{lem2} we get
\begin{corollary}\label{cor1}
 For all $j\in\nn$, all $k\in\{1,\ldots, n_j\},$ all $\alpha$ and $\beta$ in $\nn$, $l=1,2$ and all $\zeta\in\pk1\zjk$ we have uniformly with respect to $j,k,l$ and $\zeta$
$$\left|
\frac{P^{(j,k)}_l(\zeta)}{f_l(\zeta)}\diffp{^{\alpha+\beta}f_l}{{\zeta^*_1}^\alpha\partial {\zeta_2^*}^\beta}(\zeta)
\right|\leqs |\rho(\zjk)|^{i_{l}^{(j,k)}-\alpha-\frac\beta2}.$$
\end{corollary}
In the proof of the following corollary appears the technical reason why we have to introduce the open sets ${\cal U}_1^{j,k}$ and ${\cal U}_2^{j,k}$.
\begin{corollary}\label{cor3}
 For $l,m\in\{1,2\}$, we can write $\frac{P^{(j,k)}_l}{f_l} b_m=\varphi_1^{(j,k,l,m)}d\zeta^*_1+\varphi_2^{(j,k,l,m)}d\zeta^*_2$ with $\varphi_1^{(j,k,l,m)}$ and $\varphi_2^{(j,k,l,m)}$ satisfying for all $\zeta\in{\cal U}^{(j,k)}_l$ and  all differential operator $\nabla_z$ of order $1$ acting on $z$,
\begin{eqnarray*}
 \left|\varphi_1^{(j,k,l,m)}(\zeta,z)\right|&\leqs&\sum_{0\leq \alpha+\beta\leq \max(p_1,p_2)} \rzjk^{i_l^{(j,k)}-1}\left|\frac{\delta(\zeta,z)}{\rho(\zjk)}\right|^{\alpha+\frac\beta2},\\
\left|\varphi_2^{(j,k,l,m)}(\zeta,z)\right|&\leqs&\sum_{0\leq \alpha+\beta\leq \max(p_1,p_2)} \rzjk^{i_l^{(j,k)}-\frac12}\left|\frac{\delta(\zeta,z)}{\rho(\zjk)}\right|^{\alpha+\frac\beta2},\\
\left|\nabla_z\varphi_1^{(j,k,l,m)}(\zeta,z)\right|&\leqs&\sum_{0\leq \alpha+\beta\leq \max(p_1,p_2)} \rzjk^{i_l^{(j,k)}-2}\left|\frac{\delta(\zeta,z)}{\rho(\zjk)}\right|^{\alpha+\frac\beta2},\\
\left|\nabla_z\varphi_2^{(j,k,l,m)}(\zeta,z)\right|&\leqs&\sum_{0\leq \alpha+\beta\leq \max(p_1,p_2)} \rzjk^{i_l^{(j,k)}-\frac32}\left|\frac{\delta(\zeta,z)}{\rho(\zjk)}\right|^{\alpha+\frac\beta2},
\end{eqnarray*}
uniformly with respect to $\zeta,z,j$ and $k$.
\end{corollary}
\pr Without restriction we assume $l=1$ and for $m=1,2$, we write $b_m(\zeta,z)=b^*_{m,1}(\zeta,z)d\zeta^*_1+b^*_{m,2}(\zeta,z)d\zeta^*_2$ where $b_{m,n}^*=\int_0^1\diffp{f_m}{\zeta^*_n}(\zeta+t(z-\zeta))dt$. So
\begin{eqnarray*}
\lefteqn{b_{m,n}^*(\zeta,z)}\\
&=&\hskip-4pt\sum_{0\leq \alpha+\beta\leq\max(p_1,p_2)}\hskip-2pt \frac1{\alpha+\beta+1} \diffp{^{\alpha+\beta+1} f_m}{\zeta^*_n\partial{\zeta_1^*}^\alpha\partial{\zeta^*_2}^\beta}(\zeta)(z^*_1-\zeta^*_1)^\alpha(z^*_2-\zeta^*_2)^\beta\hskip-2pt +\hskip-1pt o\left(\hskip-1pt|z-\zeta|^{\max(p_1,p_2)}\hskip-1pt\right)
\end{eqnarray*}
and Corollary \ref{cor1} yields for all $\zeta\in\pk1\zjk$~:
$$\left| \frac{\overline{P_1^{(j,k)}(\zeta)}}{f_1(\zeta)} b_{1,1}(\zeta,z)\right|\leqs \sum_{0\leq \alpha+\beta\leq \max(p_1,p_2)} \rzjk^{i_l^{(j,k)}-1}\left|\frac{\delta(\zeta,z)}{\rho(\zjk)}\right|^{\alpha+\frac\beta2}$$
uniformly with respect to $z,\zeta,j$ and $k$. The proof of the inequality for $\left| \frac{\overline{P_1^{(j,k)}(\zeta)}}{f_1(\zeta)} b_{1,2}(\zeta,z)\right|$ is exactly the same. The one for  $\left| \frac{\overline{P_1^{(j,k)}(\zeta)}}{f_1(\zeta)} b_{2,1}(\zeta,z)\right|$ uses the definition of ${\cal U}^{(j,k)}_1$.\\
On ${\cal U}_1^{(j,k)}$, we have $\left|\frac{P_1^{(j,k)}}{f_1}\right|\leqs \left|\frac{P_2^{(j,k)}}{f_2} \right|\rzjk^{i_1^{(j,k)}-i_2^{(j,k)}}$ and again Corollary \ref{cor1} yields 
\begin{eqnarray*}
\left| \frac{\overline{P_1^{(j,k)}(\zeta)}}{f_1(\zeta)} b_{2,1}(\zeta,z)\right|&\leqs&
\left|\frac{P_2^{(j,k)}(\zeta)}{f_2(\zeta)}b_{2,1}(\zeta,z) \right|\rzjk^{i_1^{(j,k)}-i_2^{(j,k)}}\\
&\leqs& \sum_{0\leq \alpha+\beta\leq \max(p_1,p_2)} \rzjk^{i_l^{(j,k)}-1}\left|\frac{\delta(\zeta,z)}{\rho(\zjk)}\right|^{\alpha+\frac\beta2} 
\end{eqnarray*}
uniformly with respect to $z,\zeta,j$ and $k$. Again, the inequality for $\left| \frac{\overline{P_1^{(j,k)}(\zeta)}}{f_1(\zeta)} b_{2,2}(\zeta,z)\right|$ can be obtained in the same way.\qed
\par\medskip
Corollary \ref{cor3} and \ref{cor4} imply  for some $N'$ arbitrarily large provided $N$ is large enough, that

\begin{eqnarray*}
 \lefteqn{\left|
 \frac{\overline{P_l^{(j,k)}}(\zeta)}{f_l(\zeta)} b_m(\zeta,z)\wedge \diffp{^{i_l^{(j,k)}}}{\overline{\zeta^*_2}^{i^{(j,k)}_l}} \left( \chi_l^{(j,k)}(\zeta) \overline\partial h(\zeta)P^{N,1}(\zeta,z)\right)
\right|}\\&&\hskip150pt \leq
 \rzjk^{-3}\left(\frac{\rzjk}{\rzjk+|\rho(z)|+\delta(\zjk,z)}\right)^{N'} \tilde h(\zeta)
\end{eqnarray*}
and for $\nabla_z$ a differential of order $1$
\begin{eqnarray*}
 \lefteqn{\left|\nabla_z\left(
 \frac{\overline{P_l^{(j,k)}}(\zeta)}{f_l(\zeta)} b_m(\zeta,z)\wedge \diffp{^{i_l^{(j,k)}}}{\overline{\zeta^*_2}^{i^{(j,k)}_l}} \left( \chi_l^{(j,k)}(\zeta) \overline\partial h(\zeta)P^{N,1}(\zeta,z)\right)\right)
\right|}\\&&\hskip125pt \leq
 |\rho(z)|^{-1} \rzjk^{-3}\left(\frac{\rzjk}{\rzjk+|\rho(z)|+\delta(\zjk,z)}\right)^{N'} \tilde h(\zeta)
\end{eqnarray*}
where $\tilde h(\zeta)=\max_{n\in\{0,\ldots, i_l^{(j,k)}\}} \left(\left|\diffp{^{n+1} h}{\overline{\zeta^*_2}^{n+1}}(\zeta)|\rho(\zeta)|^{\frac {n+1} 2} \right|, \left|\diffp{^{n+1} h}{\overline{\zeta^*_1}\partial \overline{\zeta^*_2}^{n}}(\zeta)|\rho(\zeta)|^{\frac {n} 2+1} \right|\right)$. We conclude as in \cite{Al-Maz} that Theorem \ref{main_result} holds true.

\section{Local division}\label{section6}
\subsection{Local holomorphic division}
In this subsection we will prove Theorem \ref{th3} and his analogue in the $L^q$ case, the following
\begin{theorem}\label{th4}
 When $n=2$, let $g$ be a holomorphic function defined on $D$. Assume that $X_1\cap X_2$ is a complete intersection and that there exist $\kappa>0$, a real number $q\geq 1$ and a locally finite covering $\left({\cal P}_{\kappa |\rho(\zeta_j)|}(\zeta_j)\right)_{j\in I}$ of $D$ such that for all $j\in I$, 
  there exist two function $\hat g^{(j)}_1$ and $\hat g^{(j)}_2$, $C^\infty$-smooth on ${\cal P}_{\kappa |\rho(\zeta_j)|}(\zeta_j)$, such that
  \begin{enumerate}[(a)]
   \item $g=\hat g^{(j)}_1f_1+\hat g^{(j)}_2 f_2$ on ${\cal P}_{\kappa |\rho(\zeta_j)|}(\zeta_j)$;
   \item  \label{lp2}
   $c_{l,\alpha,\beta}:=\sum_{j\in I} \int_{{\cal P}_{\kappa |\rho(\zeta_j)|}(\zeta_j)} \left|\diffp{^{\alpha+\beta}\hat g^{(j)}_l}{\overline{\zeta^*_1}^\alpha\partial\overline{\zeta^*_2}^\beta}(z)\right|^q|\rho(z_j)|^{\alpha+\frac\beta2}dV(z)<\infty$ for $l=1$ and $l=2$ and all integers $\alpha$ and $\beta$;
  \item \label{lp3} for $l=1$ and $l=2$, for all non negatives integers $\alpha,\overline\alpha,\beta$ and $\overline\beta$, there exist $N\in\nn$
 and $c>0$
 such that for all $j$, $\sup_{\pk1z}\left|\diffp{^{\alpha+\overline\alpha+\beta+\overline\beta} \hat g^{(j)}_l} {{\zeta^*_1}^\alpha\partial{\zeta^*_2}^\beta\partial{\zeta^*_1}^{\overline\alpha}\partial{\zeta^*_2}^{\overline\beta}}\right| |\rho(z)|^N|\leq c$.
  \end{enumerate}
 Then there exist two smooth functions $\tilde g_1$ and $\tilde g_2$  which satisfy {\it (\ref{mth1})-(\ref{mth3})} of Theorem \ref{main_result} with $q$.
\end{theorem}
\pr
It suffices to glue together all the $\hat g_{1}^{(j)}$ and $\hat g_{2}^{(j)}$ using a suitable partition of unity.  Let $(\chi_j)_{j\in\nn}$ be a partition of unity subordinated to
$\left({\cal P}_{\kappa |\rho(\zeta_j)|}(\zeta_i)\right)_{j\in\nn}$ such that 
for all $j$ and all
$\zeta\in{\cal P}_{\kappa |\rho(\zeta_j)|}(\zeta_j)$, we have
$\left|\diffp{^{\alpha+\overline\alpha+\beta+\overline\beta} \chi_j}
{{z^{*}_1}^\alpha\partial {z^{*}_2}^\beta\partial
\overline{z^{*}_1 }^{\overline\alpha}\partial
\overline{z^{*}_2 }^{\overline\beta}}
(\zeta)\right|\leqs \frac{1}{|\rho(\zeta_j)|^{\alpha+\overline\alpha+\frac{\beta+\overline\beta}{2}}}$, uniformly with respect to $\zeta_j$ and $\zeta$.
We set $\tilde g_1=\sum_j \chi_j \hat g^{(j)}_{1}$ and $\tilde g_2=\sum_j \chi_j \hat g^{(j)}_{2}$ and thus we get the two functions defined on $D$ which satisfy {\it (\ref{mth1})}, {\it (\ref{mth2})} and {\it (\ref{mth3})} by construction.\qed
\par\medskip
Theorem \ref{th3} can be proved in exactly the same way than Theorem \ref{th4}, so we omit the proof.
\subsection{Divided differences and division}
In order to apply Theorem \ref{th3} and \ref{th4}, we will use divided differences and find numerical conditions on $g$ which ensure the existence of local smooth division formula in $L^\infty$ or in $L^q$. We define the divided differences using the following settings.

We set
$$\Lambda^{(1)}_{z,v}=\{\lambda\in\cc,\ |\lambda|<\tau(z,v,3\kappa|\rho(z)|)\text{ and } z+\lambda v\in X_2\setminus X_1\}$$
The points $z+\lambda v$, $\lambda\in\Lambda^{(1)}_{z,v},$ are the points of $X_2\setminus X_1$ which belong to $\Delta_{z,v}\left( \tau(z,v,3\kappa|\rho(z)|)\right)$, thus they all belong to $D$ as soon as $\kappa<\frac13$. We analogously define
$$\Lambda^{(2)}_{z,v}=\{\lambda\in\cc,\ |\lambda|<\tau(z,v,3\kappa|\rho(z)|)\text{ and } z+\lambda v\in X_1\setminus X_2\}.$$
For a function $h$ defined on a subset $\cal U$ of $\cc^n$, $z\in\cc^n$, $v$ a unit vector of $\cc^n$ and $\lambda\in\cc$ such that $z+\lambda v$ belongs to $\cal U$, we set $h_{z,v}[\lambda]=h(z+\lambda v).$ If $h_{z,v}[\lambda_1,\ldots, \lambda_k]$ is defined, for 
$\lambda_1,\ldots, \lambda_{k+1}\in\cc$ pairwise distinct such that $z+\lambda_i v$ belongs to ${\cal U}$ for all $i$, we set
$$h_{z,v}[\lambda_1,\ldots,\lambda_{k+1}]:=\frac{h_{z,v}[\lambda_1,\ldots,\lambda_{k}]-h_{z,v}[\lambda_2,\ldots,\lambda_{k+1}]}{\lambda_1-\lambda_{k+1}}.$$
Now, for $z\in X_2\setminus X_1$ (resp. $z\in X_1\setminus X_2$) let us define $g^{(2)}(z)=\frac{g(z)}{f_2(z)}$ (resp. $g^{(1)}(z)=\frac{g(z)}{f_1(z)}$). For $l=1$ or $l=2$, the quantity 
$g^{(l)}_{z,v}[\lambda_1,\ldots,\lambda_k]$  make sense for all $\lambda_1,\ldots,\lambda_k\in \Lambda_{z,v}^{(l)}$ pairwise distinct.

We first prove a technical result we will need in this section.
\begin{lemma}\label{prod_dif}
 Let $\alpha$ and $\beta$ be two functions defined on a subset $\cal U$ of $\cc$. Then, for all $z_1,\ldots, z_n$ pairwise distinct points of $\cal U$ we have
$$(\alpha\cdot \beta)[z_1,\ldots, z_n]=\sum_{k=1}^n\alpha[z_1,\ldots, z_k]\cdot \beta[z_k,\cdots,z_n].$$
\end{lemma}
\pr We prove the lemma by induction on $n$, the case $n=1$ being trivial. We assume the lemma proved for $n$ points, $n\geq 1$. Let $z_1,\ldots, z_{n+1}$ be $n+1$ points of $\cal U$. Then
\begin{eqnarray*}
 \lefteqn{(\alpha\cdot\beta)[z_1,\ldots, z_{n+1}]}\\
&=&\frac{(\alpha\cdot\beta)[z_1,z_3,\ldots, z_{n+1}]-(\alpha\cdot\beta)[z_2,\ldots, z_{n+1}]}{z_1-z_2}\\
&=&\frac1{z_1-z_2} \left(\sum_{k=3}^{n+1} \alpha[z_1,z_3,\ldots z_{k}] \beta[z_k,\ldots, z_{n+1}]+\alpha[z_1]\beta[z_3,\ldots,z_{n+1}]\right)\\
&& -\frac1{z_1-z_2} \sum_{k=2}^{n+1} \alpha[z_2,\ldots z_{k}] \beta[z_k,\ldots, z_{n+1}]\\
&=&\sum_{k=3}^{n+1} \frac{\alpha[z_1,z_3,\ldots z_{k}]-\alpha[z_2,\ldots z_{k}]}{z_1-z_2} \beta[z_k,\ldots, z_{n+1}]+\\
&&\frac{\alpha[z_1]-\alpha[z_2]}{z_1-z_2} \beta[z_2,\ldots,z_{n+1}] + 
\alpha[z_1]\frac{\beta[z_1,z_3,\ldots,z_{n+1}]- \beta[z_2,\ldots,z_{n+1}]}{z_1-z_2}.\qed
\end{eqnarray*}

\subsubsection{The $L^\infty-BMO$-case}\label{BMOcase}
In this subsection, we establish the necessary conditions in $\cc^n$ and the sufficient conditions $\cc^2$ for a function $g$ to be written as $g=g_1f_1+g_2f_2$, $g_1$ and $g_2$ smooth functions satisfying the hypothesis of Theorem \ref{main_result}.\\
For $l=1$ and $l=2$ let us define the numbers
$$c^{(l)}_\infty(g)=\sup\left(|g^{(l)}_{z,v}[\lambda_1,\ldots,\lambda_k]|  \tau(z,v,|\rho(z)|)^{k-1}\right)$$
where the supremum is taken over all $z\in D,$ all $v\in\cc^n$ with $|v|=1$ and all $\lambda_1,\ldots,\lambda_k\in\Lambda^{(l)}_{z,v}$ pairwise distinct.\\
We have the following necessary conditions in $\cc^n$, $n\geq 2$.
\begin{theorem}\label{th1}
In $\cc^n$, $n\geq 2$, let $g$ be a holomorphic on $D$ and let $g_1,g_2$ be two bounded holomorphic functions on $D$ such that $g=g_1f_1+g_2f_2$. Then
$$\left\| \frac{g}{\max(|f_1|,|f_2|)}\right\|_{L^\infty(D)} \leqs \max(\|g_1\|_{L^\infty(D)},\|g_2\|_{L^\infty(D)})$$
and for $l=1,2$ :
$$c^{(l)}_\infty(g) \leqs \sup_{b\Delta_{z,v} \left(4\kappa\tau(z,v,|\rho(z)|)\right)}|g_l|.$$
\end{theorem}
\noindent{\it Proof~:}
The first point is trivial and we only prove the second one for $l=1$. Let $\lambda_1,\ldots,\lambda_k$ be $k$ pairwise distinct elements of $\Lambda^{(1)}_{z,v}$. For all $i$ we have $g^{(1)}_{z,v}[\lambda_i]=g_1(z+\lambda_i v)$ because $f_2(z+\lambda_i v)=0$. Therefore, $g_{z,v}^{(1)}[\lambda_1,\ldots, \lambda_k]={(g_1)}_{z,v}[\lambda_1,\ldots, \lambda_k]$. As in \cite{Al-Maz}, it then follows from Cauchy's formula
that
\begin{eqnarray*}
 |g^{(1)}_{z,v}[\lambda_1,\ldots,\lambda_k]|&\leqs& \left| \frac1{2i\pi}\int_{|\lambda|= \tau(z,v,4\kappa|\rho(z)|)} \frac{g_1(z+\lambda v)}{\prod_{i=1}^k(\xi -\lambda_i)}d\xi\right|\\
&\leqs&  \tau(z,v,|\rho(z)|)^{-k+1}\sup_{b\Delta_{z,v} \left(4\kappa\tau(z,v,|\rho(z)|)\right)}|g_1|.\qed
\end{eqnarray*}
Now we prove that these conditions are sufficient in $\cc^2$ in order to get a $BMO$ division.
\begin{theorem}\label{th2}
In $\cc^2$, let $g$ be a holomorphic function on $D$ which belong to the ideal generated by $f_1$ and $f_2$ and such that
\begin{enumerate}[(i)]
 \item \label{hyp1} $c(g)=\sup_{z\in D} \frac{|g(z)|}{\max(|f_1(z)|,|f_2(z)|)}<\infty$,
 \item \label{hyp2} $c^{(1)}_\infty(g)$ and $c^{(2)}_\infty(g)$ are finished.
\end{enumerate}
Then for all $z\in D$, there exist two holomorphic functions $g_1$ and $g_2$ which beblong to $BMO(D)$ and such that $g_1f_1+g_2f_2=g$.
\end{theorem}
\noindent{\it Proof~:}
It suffices to construct for all $z$ near $bD$ two smooth functions $\hat g_1$ and $\hat g_2$ on ${\cal P}_{\kappa|\rho(z)|}(z)$ which satisfy $(\ref{th3i})$ and $(\ref{th3ii})$ of Theorem \ref{th3}. 

Let $\zeta_0$ be a point in $bD$. If $f_1(\zeta_0)\neq 0$ then $f_1$ does not vanish on a neighborhood ${\cal U}_0$ of $\zeta_0$. Then we can define $\hat g_1= \frac g{f_1}$, $\hat g_2=0$ which obviously satisfy $(\ref{th3i})$ and $(\ref{th3ii})$ for all $z\in D$ close to  $\zeta_0$. We proceed analogously if $f_2(\zeta_0)\neq 0$.

If $\zeta_0$ belongs to $X_1\cap X_2\cap bD$, since the intersection $X_1\cap X_2$ is complete, without restriction we can choose a neighborhood ${\cal U}_0$ of $\zeta_0$ such that $X_1\cap X_2\cap {\cal U}_0=\{\zeta_0\}$. Then we fix some point $z$ in ${\cal U}_0$ and we construct $\hat g_1$ and $\hat g_2$ on ${\cal P}_{\kappa|\rho(z)|}(z)$ which satisfy $(\ref{th3i})$ and $(\ref{th3ii})$ of Theorem \ref{th3}. We denote by $p_1$ and $p_2$ the order of $\zeta_0$ as zero of $f_1$ and $f_2$ respectively.
We also denote by $(\zeta^*_{0,1},\zeta^*_{0,2})$ the coordinates of $\zeta_0^*$ in the \ko coordinates at $z$. If $|\zeta_{0,1}^*|<2\kappa |\rho (z)|$, then for $l=1$ and $l=2$ we set $i_l=0$, $P_l(\zeta)=1$ and $Q_l(\zeta)=f_l(\zeta)$. Otherwise, we use the parametrization 
  $\alpha_{1,i}$, $i\in\{1,\ldots,p_1\}$, of $X_1$ and $\alpha_{2,i}$, $i\in\{1,\ldots,p_2\}$, of $X_2$
   given by Proposition \ref{propII.0.2}. We denote by $I_{l}$ the set $I_{l}=\{i, \exists z^*_1\in \Delta_0(\kappa |\rho(z)|)\text{ such that } |\alpha_{l,i}(z^*_1)|\leq (2\kappa |\rho(z)|)^{\frac12}\}$, $i_{l}=\# I_{l}$, $P_{l}(\zeta)=\prod_{i\in I_{l}} (\zeta_2^*-\alpha_{l,i}(\zeta^*_1))$ and $Q_{l}(\zeta)=\frac {f_l}{P_{l}}$.

If $i_1=0$ we set $\tilde g_{2}=0$. Otherwise, without restriction we assume that $I_{1}=\{1,\ldots, i_{1}\}$ and for $k\leq i_{1}$ and $\zeta^*_1$ such that $f_2(z+\zeta^*_1\eta_{z} +\alpha_{1,i}(\zeta^*_1)v_{z})\neq 0$, we introduce $g^{(2)}=\frac g{f_2}$ and
\begin{eqnarray}
\check g_{1,\ldots,k}^{(2)} (\zeta^*_1):=\left(\frac g{P_2}\right)_{z+\zeta^*_1\eta_z,v_z}[\alpha_{1,1}(\zeta^*_1),\ldots,\alpha_{1,k}(\zeta^*_1)] .\label{eq20}
\end{eqnarray}
Since $X_1\cap X_2\cap {\cal U}_0=\{\zeta_0\}$, $\check g_{1,\ldots,k}^{(2)}$ is defined on $\Delta_0(2\kappa |\rho(z)|)$ and we have by Lemma \ref{prod_dif}
\begin{eqnarray*}
 \lefteqn{\check g_{1,\ldots,k}^{(2)} (\zeta^*_1)}\\
&& = \left( \frac g{P_{2}} \right)_{z+\zeta^*_1 \eta_{z}, v_{z}}[\alpha_{1,1}(\zeta^*_1),\ldots, \alpha_{1,k}(\zeta^*_1)]                                    \\
&& = \left( \frac g {f_2} Q_{2} \right)_{z+\zeta^*_1 \eta_{z}, v_{z}}[\alpha_{1,1}(\zeta^*_1),\ldots, \alpha_{1,k}(\zeta^*_1)]\\
&& = \sum_{j=1}^k g^{(2)}_{z+\zeta^*_1 \eta_{z}, v_{z}} [\alpha_{1,1}(\zeta^*_1),\ldots, \alpha_{1,j}(\zeta^*_1)]
\left(Q_{2} \right)_{z+\zeta^*_1 \eta_{z}, v_{z}} [\alpha_{1,j}(\zeta^*_1),\ldots, \alpha_{1,k}(\zeta^*_1)].
\end{eqnarray*}
Now from \cite{Mon} we have 
$$|\left(Q_{2} \right)_{z+\zeta^*_1 \eta_{z}, v_{z}} [\alpha_{1,j}(\zeta^*_1),\ldots, \alpha_{1,k}(\zeta^*_1)]|\leqs |\rho(z)|^{\frac{j-k}2} \sup_{|\xi|=(4\kappa |\rho(z)|)^{\frac12}} |Q_2(z+\zeta^*_1\eta_z+\xi v_z)|,$$ 
which, with the assumption $c^{(2)}_\infty(g)<\infty$,  gives for all $\zeta^*_1\in  \Delta_0(2\kappa |\rho(z)|)$~:
\begin{eqnarray}
|{\check g_{1,\ldots,k}^{(2)} (\zeta^*_1)}|
&\leqs& c^{(2)}_\infty(g) |\rho(z)|^{\frac{1-k}2}\sup_{|\xi|=(4\kappa |\rho(z)|)^{\frac12}} |Q_2(z+\zeta^*_1\eta_z+\xi v_z)|.\label{eq10}
\end{eqnarray}
Now we set 
$$\tilde g_{2}(\zeta)=\sum_{k=1}^{i_{2}} \check g_{1,\ldots, k}^{(2)}(\zeta^*_1) \prod_{i=1}^{k-1} (\zeta^*_2-\alpha_{1,i}(\zeta^*_1)).$$
and we define $\tilde g_{1}$ analogously. For $\zeta_1^*$ fixed, $\tilde g_2(\zeta^*_1,\cdot)$ is the polynomial  which interpolates $\frac{g(\zeta^*_1,\cdot)}{P_2(\zeta^*_1,\cdot)}$ at the points $\zeta^*_2=\alpha_{1,1}(\zeta^*_1),\ldots, \alpha_{1,i_1}(\zeta_1^*)$.\\
Since $|\alpha_{1,i}(\zeta^*_1)|\leqs |\rho(z)|^{\frac12}$ for all $i\in I_1$ and all $\zeta^*_1\in\Delta_0(2\kappa |\rho(z)|)$, (\ref{eq10}) yields for all $\zeta\in{\cal P}_{2\kappa|\rho(z)|}({z})$~:
\begin{eqnarray}
|\tilde g_{2}(\zeta)|\leqs c^{(2)}_\infty(g)\sup_{\over{|\xi_2|\leq(4\kappa |\rho(z)|)^{\frac12}}{|\xi_1|\leq2\kappa |\rho(z)|}} |Q_2(z+\xi_1\eta_z+\xi_2 v_z)|.\label{eq11}
\end{eqnarray}
Then Cauchy's inequalities gives for all $\zeta\in{\cal P}_{\kappa|\rho(z)|}({z})$ and all $\alpha$ and $\beta$
\begin{eqnarray}
\left|\diffp{^{\alpha+\beta}\tilde g_{2}}{{\zeta^*_1}^\alpha\partial{\zeta^*_2}^\beta}(\zeta)\right|\leqs c^{(2)}_\infty(g)|\rho(z)|^{-\alpha-\frac\beta2}\sup_{\over{|\xi_2|\leq(4\kappa |\rho(z)|)^{\frac12}}{|\xi_1|\leq2\kappa |\rho(z)|}} |Q_2(z+\xi_1\eta_z+\xi_2 v_z)|.\label{eq11bis}
\end{eqnarray}

Now let $g_1$ and $g_2$ be holomorphic functions on $D$ such that $g=f_1g_1+f_2g_2$. Then $\tilde g_{2}(\zeta^*_1,\cdot)$ interpolates $g_{2} (\zeta_1^*,\cdot) Q_{2}(\zeta_1^*,\cdot)$ at the points $\alpha_{1,i}(\zeta^*_1)$ for all $i\in I_{2}$ because for such an $i$ we have by definition
\begin{eqnarray*}
\tilde g_{2}(\zeta^*_1,\alpha_{1,i}(\zeta^*_1))&=&\frac{g(\zeta^*_1,\alpha_{1,i}(\zeta^*_1)}{P_{2}(\zeta^*_1,\alpha_{1,i}(\zeta^*_1))}\\
&=& g_2(\zeta^*_1,\alpha_{1,i}(\zeta^*_1))\cdot Q_{2}(\zeta^*_1,\alpha_{1,i}(\zeta^*_1))\\
\end{eqnarray*}
Therefore we can write
\begin{eqnarray}
g_2(\zeta)=\frac1{Q_{2}(\zeta)}\left(\tilde g_{2}(\zeta)+P_{1}(\zeta) \cdot e_1(\zeta)\right) \label{eq22} 
\end{eqnarray}
where $e_1$ is the interpolation error which is given by
\begin{eqnarray}
e_1(\zeta)=\frac1{2i\pi} \int_{|\xi|=(4\kappa|\rho(z)|)^{\frac12}} \frac{g_2(\zeta_1^*,\xi) Q_{2} (\zeta_1^*,\xi)}{P_{1}(\zeta_1^*,\xi)\cdot (\xi-\zeta^*_2)} d\xi.\label{eq23}
\end{eqnarray}
We have an analogous expression for $g_1$. We point out that (\ref{eq22}) and its analogous for $g_1$ also holds if $i_1=0$ or $i_2=0$.\\
This yields
\begin{eqnarray}
g(\zeta)&=&f_1(\zeta)g_1(\zeta) +f_2(\zeta)g_2(\zeta)\nonumber\\
&=& P_{1}(\zeta)\tilde g_{1}(\zeta)+P_{2}(\zeta)\tilde g_{2}(\zeta) +P_{1}(\zeta)P_{2}(\zeta) e(\zeta)\label{eq24}
\end{eqnarray}
where
\begin{eqnarray*}
e(\zeta)&=&e_1(\zeta)+e_2(\zeta)\nonumber\\
&=& \frac1{2i\pi}\int_{|\xi|=(4\kappa|\rho(z)|)^{\frac12}} \frac{g(\zeta_1^*,\xi)}{P_{1}(\zeta_1^*,\xi)\cdot P_{2}(\zeta_1^*,\xi)\cdot (\xi-\zeta^*_2)} d\xi.
 \end{eqnarray*}
Searching for $\hat g_1$ and $\hat g_2$ such that $g=\hat g_1f_1+\hat g_2f_2$ in $\pk1z$, since there is a factor $P_1P_2$ in front $e$ in (\ref{eq24}), we can put $P_2e$ either in $\hat g_1$ with $\tilde g_1$ or we can put $P_1e$ in $\hat g_2$ with $\tilde g_2$. But in order to have a good upper bound, we have to cut it in to two pieces in a suitable way. This will be done analogously to the construction of the currents. Let
\begin{eqnarray*}
 {\cal U}_{1}&:=&\left\{\zeta\in \pk1z,\ \left|\frac{f_1(\zeta)\rho(z)^{i_1}}{{P}_1(\zeta)}\right|>\frac13\left|\frac{f_2(\zeta)\rho(z)^{i_2}}{P_2(\zeta)}\right|\right\},\\
 {\cal U}_{2}&:=&\left\{\zeta\in \pk1z,\frac23 \left|\frac{f_2(\zeta)\rho(z)^{i_2}}{P_2(\zeta)}\right|>\left|\frac{f_1(\zeta)\rho(z)^{i_1}}{{P}_1(\zeta)}\right|\right\}.
\end{eqnarray*}
Let also $\chi $ be a smooth function on $\cc^2\setminus\{0\}$  such that $\chi (z_1,z_2)=1$ if $|z_1|>\frac23 |z_2|$ and $\chi (z_1,z_2)=0$ if $|z_1|<\frac13|z_2|$. We set $\chi_1(\zeta)=\chi\left(\frac{f_1(\zeta)\rho(z)^{i_1}}{P_1(\zeta)},\frac{f_2(\zeta)\rho(z)^{i_2}}{P_2(\zeta)}\right)$, $\chi_2(\zeta)=1-\chi_1(\zeta)$ and 
\begin{eqnarray*}
{\hat g_{1}(\zeta)}
&=&\frac1{Q_{1}(\zeta)}\left(\tilde g_{1}(\zeta)+\chi_1(\zeta){P_{2}(\zeta)}e(\zeta)\right),\\
{\hat g_{2}(\zeta)}
&=&\frac1{Q_{2}(\zeta)}\left(\tilde g_{2}(\zeta)+\chi_2(\zeta){P_{1}(\zeta)}e(\zeta)\right).
\end{eqnarray*}
Since $Q_2=\frac{f_2}{P_2}$, Lemma \ref{lem1} and (\ref{eq11bis}) give for all $\zeta\in{\cal P}_{\kappa|\rho(z)|}({z})$~:
\begin{eqnarray}
\left|\diffp{^{\alpha+\beta}}{{\zeta^*_1}^\alpha \partial{\zeta^*_2}^\beta}\left(\frac1{Q_{2}(\zeta)}\tilde g_{2}(\zeta)\right)\right| &\leqs& c_\infty^{(2)}(g) |\rho(z)|^{-\alpha-\frac\beta2}\label{eq22bis}
\end{eqnarray}
From assumption (\ref{hyp1}), We get for all $\zeta_1^*\in \Delta_0(2\kappa|\rho(z)|)$ and all $\xi$ such that $|\xi|=(4\kappa|\rho(z)|)^{\frac12}$~:
\begin{eqnarray*}
 \left|{g(\zeta^*_1,\xi)}\right|&\leq& c(g)\left(\sup_{\over{|\xi^*_1|\leq 2\kappa|\rho(z)|}{|\xi^*_2|\leq (4\kappa|\rho(z)|)^\frac12}} |f_1(\xi)|+ \sup_{\over{|\xi^*_1|\leq 2\kappa|\rho(z)|}{|\xi^*_2|\leq (4\kappa|\rho(z)|)^\frac12}} |f_2(\xi)| \right).
\end{eqnarray*}
And so Cauchy's inequalities yield for all integer $\alpha$ and all $\zeta_1^*\in \Delta_0(\kappa|\rho(z)|)$ and all $\xi$ such that $|\xi|=(4\kappa|\rho(z)|)^{\frac12}$
\begin{eqnarray*}
 \lefteqn{\left|{\diffp{^{\alpha} g}{{\zeta^*_1}^\alpha}(\zeta^*_1,\xi)}\right|}\\
 &\leqs& c(g)|\rho(z)|^{-\alpha}
\left(|\rho(z)|^{i_1}\sup_{\over{|\xi^*_1|\leq 2\kappa|\rho(z)|}{|\xi^*_2|\leq (4\kappa|\rho(z)|)^\frac12}} |Q_1(\xi)|+ |\rho(z)|^{i_2}\sup_{\over{|\xi^*_1|\leq 2\kappa|\rho(z)|}{|\xi^*_2|\leq (4\kappa|\rho(z)|)^\frac12}} |Q_2(\xi)| \right).
\end{eqnarray*}
Therefore, for all $\zeta\in  {\cal U}_{2}$ and all non negative integers $\alpha$ and $\beta$, Lemma \ref{lem1} gives 
\begin{eqnarray*}
 \left|\diffp{^{\alpha+\beta}}{{\zeta^*_1}^\alpha \partial{\zeta^*_2}^\beta} \left(\frac{P_{1}(\zeta)}{Q_2(\zeta)}  e(\zeta)\right) \right|&\leqs&|\rho(z)|^{-\alpha-\frac\beta2}c(g).
\end{eqnarray*}
Since for all $\alpha,\beta \in\nn$, $\left|\diffp{^{\alpha+\beta}\chi_2}{\overline{\zeta^*_1}^\alpha\partial \overline{\zeta^*_2}^\beta} (z) \right|\leqs |\rho(z_j)|^{-\alpha-\frac\beta2}$ (see Lemma \ref{lem1}), with (\ref{eq22bis}), this yields
\begin{eqnarray*}
{\left|\diffp{^{\alpha+\overline\alpha+\beta+\overline\beta} \hat g_{2}}{{\zeta^*_1}^\alpha\partial{\zeta^*_2}^\beta\partial \overline{\zeta^*_1}^{\overline\alpha}\partial\overline{\zeta^*_2}^{\overline\beta}} (\zeta)\right|}
&\leqs& 
|\rho(z)|^{-\alpha-\overline\alpha-\frac{\beta+\overline\beta}2}\left(c(g)+c^{(2)}_\infty(g)\right).
\end{eqnarray*}
The same inequality holds for $\hat g_1$ and we have finally proved that $\hat g_1$ and $\hat g_2$ are smooth functions such that 
(\ref{th3i}) and (\ref{th3ii}) of Theorem \ref{th3} hold.\qed 

\subsection{The $L^q$-case}
The assumption under which a function $g$ holomorphic on $D$ can be written as $g=g_1f_1+g_2f_2$ with $g_1$ and $g_2$ being holomorphic on $D$ and belonging to $L^q(D)$ uses a $\kappa $-covering $\left({\cal P}_{\kappa |\rho(z_j)|}(z_j)\right)_{j\in\nn}$ in addition to the divided differences.\\
By transversality of $X_1$ and $bD$, and of $X_2$ and $bD$, for all $j$ there exists $w_j$ in the complex tangent plane to $bD_{\rho(z_j)}$ such that $\pi_j$, the orthogonal projection on the hyperplane orthogonal to $w_j$ passing through $z_j$, is a covering of $X_1$ and $X_2$. We denote by $w_1^*,\ldots, w^*_n$ an orthonormal basis of $\cc^n$ such that $w_1^*=\eta_{z_j}$ and $w_n^*=w_j$ and we set ${\cal P}'_{\varepsilon}(z_j)=\{z'=z_j+z^*_1w^*_1+\ldots+z^*_{n-1} w^*_{n-1},\ |z^*_1|< \varepsilon \text{ and } |z_k^*|<\varepsilon^{\frac12},\ k=2,\ldots, n-1\}$. 
We put
\begin{eqnarray*}
c^{(l)}_{q,\kappa,{(z_j)_{j\in\nn}}}(g)\hskip-5pt &=&\hskip-5pt\sum_{j=0}^\infty \int_{z'\in{\cal P}'_{2\kappa |\rho(z_j)|}(z_j)}
\sum_{\over{\lambda_1,\ldots,\lambda_k\in\Lambda_{z',w_n^*}}{\lambda_i\neq\lambda_l\text{ for }i\neq l}}\hskip - 3pt
|\rho(z_j)|^{q\frac{k-1}2+1} \left|g^{(l)}_{z',w_n^*}[\lambda_1,\ldots,\lambda_k]\right| dV_{n-1}(z')
\end{eqnarray*}
where $dV_{n-1}$ is the Lebesgue measure in $\cc^{n-1}$ and $g^{(l)}=\frac{g}{f_l}$, $l=1$ or $l=2$.

Now we prove the following necessary conditions
\begin{theorem}\label{th5}
 Let $g_1$ and $g_2$ belonging to $L^q(D)$ be two holomorphic functions on $D$ and set $g=g_1f_1+g_2f_2$. Then
\begin{enumerate}[(i)]
 \item \label{second_point} $\frac{g}{\max(|f_1|,|f_2|)}$ belongs to $L^q(D)$ and $\left\| \frac{g}{\max(|f_1|,|f_2|)}\right\|_{L^q(D)}\leqs \max(\|g_1\|_{L^q(D)},\|g_2\|_{L^q(D)})$.
\item \label{first_point} For $l=1$ or $l=2$ and any $\kappa$-covering $\left({\cal P}_{\kappa|\rho(z_j)|} (z_j)\right)_j$, we have $c_{q,\kappa,(z_j)_j}^{(l)}(g)\leqs \|g_l\|^q_{L^q(D)}$,
\end{enumerate}
\end{theorem}
\pr The point (\ref{second_point}) is trivial and we only prove (\ref{first_point}). As in the proof of Theorem \ref{th1}, for all $j\in\nn$, all $z'\in {\cal P}'_{\kappa|\rho(z_j)|}(z_j)$ and all $r\in[\frac72\kappa|\rho(z_j)|^{\frac12},4\kappa|\rho(z_j)|^{\frac12}]$ we have
\begin{eqnarray*}
 g^{(l)}_{z',w^*_n}[\lambda_1,\ldots,\lambda_k]&=& \frac1{2i\pi} \int_{|\lambda|=r} \frac{g_l(z'+\lambda w_n^*)}{\prod_{i=1}^k(\xi -\lambda_i)}d\xi.
\end{eqnarray*}
After integration for $r\in[(7/2\kappa |\rho (z_j)|)^{\frac12},(4\kappa |\rho (z_j)|)^{\frac12}]$, 
Jensen's inequality yields
$$\left|g^{(l)}_{z',w_n^*}[\lambda_1,\ldots,\lambda_k]\right|^q\leqs |\rho(z_j)|^{\frac{1-k}2 q-1} \int_{|\lambda|\leq (4\kappa |\rho(z_j)|)^\frac12} |g_l(z'+\lambda w^*_n)|^q dV_1(\lambda).$$
Now we integrate the former inequality for $z'\in{\cal P}'_{\kappa |\rho(z_j)|}(z_j)$ and get
\begin{eqnarray*}
{\int_{z'\in {\cal P}'_{\kappa|\rho(z_j)|}(z_j)} \left|g^{(l)}_{z',w^*_n}[\lambda_1,\ldots,\lambda_k]\right|^q |\rho(z_j)|^{\frac{k-1}2 q+1}dV_{n-1}}
& \leqs& \int_{z\in {\cal P}_{4\kappa|\rho(z_j)|} (z_j)} |g_l(z)|^q dV_n(z). 
\end{eqnarray*}
Since $\left({\cal P}_{\kappa |\rho(z_j)|}(z_j)\right)_{j\in\nn}$ is a $\kappa $-covering, we deduce from this inequality that 
$c^{(l)}_{q,\kappa,(z_j)_{j\in\nn}}(g)\leqs \|g_l\|^q_{L^q(D)}$.\qed
\begin{theorem}\label{th6}
Let $g$ be a holomorphic function on $D$ belonging to the ideal generated by $f_1$ and $f_2$ and such that $c_{q,\kappa,(z_j)_j}^{(l)}(g)$ is finite and such that $\frac{g}{\max(|f_1|,|f_2|)}$ belongs to $L^q(D)$.\\
Then there exist two holomorphic functions $g_1$ and $g_2$ which belong to $L^q(D)$ and such that $g=g_1f_1+g_2f_2.$
\end{theorem}
\pr We aim to apply Theorem \ref{th4}. For all $j$ in $\nn$, in order to construct on $\pk1{z_{j}}$ two functions $\hat g_1^{(j)}$ and $\hat g^{(j)}_2$ which satisfy the assumption of Theorem \ref{th4}, we proceed as in the proof of Theorem \ref{th2}. The main difficulty occurs, as in the proof of Theorem \ref{th2}, when we are near a point $\zeta_0$ which belongs to $X_1\cap X_2\cap bD$.
We denote by $(\zeta_{0,1}^*,\zeta_{0,2}^*)$ the coordinates of $\zeta_0$ in the \ko coordinates at $z_{j}$.
If $|\zeta_{0,1}^*|<2\kappa|\rho(z_{j_0})|$, we set $i_{1,j}=i_{2,j}=0$, $I_{1,j}=I_{2,j}=\emptyset$, $P_{1,j}=P_{2,j}=1$, $Q_{1,j}=f_1$ and $Q_{2,j}=f_2$.
Otherwise,  we use the parametrization   $\alpha^{(j)}_{1,i}$, $i\in\{1,\ldots,p^{(j)}_1\}$ of $X_1$ and $\alpha^{(j)}_{2,i}$, $i\in\{1,\ldots,p^{(j)}_2\}$ of $X_2$ given by Proposition \ref{propII.0.2} and for $l=1$ and $l=2$, we still denote by $I_{l,j}$ the set $I_{l,j}=\{i, \exists z^*_1\in \Delta_0(\kappa |\rho(z_{j})|)\text{ such that } |\alpha^{(j)}_{l,i}(z^*_1)|\leq 2\kappa |\rho(z_j)|)^{\frac12}\}$, $i_{l,j}=\# I_{l,j}$, $P_{l,j}(\zeta)=\prod_{i\in I_{l,j}} (\zeta_2^*-\alpha^{(j)}_{l,i}(\zeta^*_1))$ and $Q_{l,j}=\frac {f_l}{P_{l,j}}$. We define $\tilde g_{1}^{(j)}$ and $\tilde g_{2}^{(j)}$ as $\tilde  g_1$ and $\tilde  g_2$ in the proof of Theorem \ref{th2}. Instead of defining $e^{(j)}_1$ and $e^{(j)}_2$ by integrals over the set $\{|\xi|=(4\kappa |\rho(z_j)|)^{\frac12}\}$ as we defined $e_1$ and $e_2$ in the proof ot Theorem \ref{th2}, here we integrate over $\{(\frac72\kappa |\rho(z_j)|)^{\frac12}\leq |\xi|\leq (4\kappa |\rho(z_j)|)^{\frac12}\}$ and set 
$$e^{(j)}(\zeta)\hskip-2pt=\hskip-2pt\frac1{2\pi(2-\sqrt {\frac72})(\kappa |\rho(z_j)|}\hskip-2pt\int_{\{(\frac72\kappa |\rho(z)|)^{\frac12}\leq |\xi|\leq (4\kappa |\rho(z)|)^{\frac12}\}}
\hskip-4pt\frac{g(z^*_1,\xi)}{P_{1,j}(z^*_1,\xi)P_{2,j}(z^*_1,\xi) (z^*_2-\xi)} dV(\xi).$$
We therefore have for all $j$ and all $z\in\pk1{z_j}$ :
$$g(z)=\tilde g^{(j)}_1(z)P_{1,j}(z)+\tilde g^{(j)}_2(z)P_{2,j}(z)+P_{1,j}(z)P_{2,j}(z)e^{(j)}(z).$$
We split $\pk1{z_j}$ in two parts as in Theorem \ref{th2} and set 
\begin{eqnarray*}
 {\cal U}^{(j)}_{1}&:=&\left\{\zeta\in \pk1{z_j},\ \left|\frac{f_1(\zeta)\rho(z_j)^{i_{1,j}}}{{P}_{1,j}(\zeta)}\right|>\frac13\left|\frac{f_2(\zeta)\rho(z_j)^{i_{2,j}}}{P_2(\zeta)}\right|\right\},\\
 {\cal U}^{(j)}_{2}&:=&\left\{\zeta\in \pk1{z_j},\frac23 \left|\frac{f_2(\zeta)\rho(z_j)^{i_{2,j}}}{P_{2,j}(\zeta)}\right|>\left|\frac{f_1(\zeta)\rho(z_j)^{i_{1,j}}}{{P}_{1,j}(\zeta)}\right|\right\}.
\end{eqnarray*}
We still denote by $\chi $ a smooth function on $\cc^2\setminus\{0\}$  such that $\chi (z_1,z_2)=1$ if $|z_1|>\frac23 |z_2|$ and $\chi (z_1,z_2)=0$ if $|z_1|<\frac13|z_2|$; and we set $\chi^{(j)}_1(\zeta)=\chi\left(\frac{f_1(\zeta)\rho(z_j)^{i_{1,j}}}{P^{(j)}_1(\zeta)},\frac{f_2(\zeta)\rho(z_j)^{i_{2,j}}}{P^{(j)}_2(\zeta)}\right)$, $\chi^{(j)}_2(\zeta)=1-\chi^{(j)}_1(\zeta)$ and 
\begin{eqnarray*}
{\hat g^{(j)}_{1}(z)}
&=&\frac1{Q^{(j)}_{1}(z)}\left(\tilde g^{(j)}_{1}(z)+\chi^{(j)}_1(z)P_{2,j}(z) e^{(j)}(z)\right),\\
{\hat g^{(j)}_{2}(z)}
&=&\frac1{Q^{(j)}_{2}(z)}\left(\tilde g^{(j)}_{2}(z)+\chi^{(j)}_2(z)P_{1,j}(z) e^{(j)}(z)\right).
\end{eqnarray*}
Therefore $g=\hat g^{(j)}_1 f_1+\hat g^{(j)}_2f_2$ on $\pk1{z_j}$ and in order to apply Theorem \ref{th4}, the assumptions (\ref{lp2}) and (\ref{lp3}) are left to be shown.\\
From Lemma \ref{lem1}, for all $j\in\nn$ and all $z\in\pk2{z_j}$, we have
$$\left|\frac1{Q_{2,j}(z)}\tilde g^{(j)}_2(z)\right|\leqs \sum_{k=1}^{i_{2,j}}|\rho(z_j)|^{\frac{k-1}2}\left|g^{(2)}_{z_j+z^*_1\eta_{z_j},v_{z_j}}[\alpha_{1,1}(z^*_1),\ldots, \alpha_{1,k}(z^*_1)]\right|$$
uniformly with respect to $z$ and $j$.\\
Therefore 
\begin{eqnarray}
\sum_{j\in\nn}
\int_{\pk2{z_j}} \left|\frac1{Q_{2,j}(z)}\tilde g^{(j)}_2(z)\right|^q dV(z)&\leqs& c^{(l)}_{q,\kappa,(z_j)}(g)\label{eq100}
\end{eqnarray}
and in particular $\frac1{Q_{2,j}}\tilde g^{(j)}_2$ is an holomorphic function with $L^q$-norm on $\pk2{z_j}$ lower than $(c^{(2)}_{q,\kappa,(z_j)}(g))^{\frac1q}$. Thus Cauchy's inequalities imply that for all $\alpha,\beta\in\nn$ and all $\pk1{z_j}$ that
\begin{eqnarray}
\left|\diffp{^{\alpha+\beta} }{{z^*_1}^\alpha\partial{z^*_2}^\beta} \left( \frac1{Q_{2,j}} \tilde g_2^{(j)}(z)\right)\right|&\leqs&  c^{(l)}_{q,\kappa,(z_j)}(g)|\rho(z_j)|^{-\alpha-\frac\beta 2}.\label{eq101}
\end{eqnarray}
Since $\frac g{\max(|f_1|,|f_2|)}$ belongs to $L^q(D)$, $g$ itself belongs to $L^q(D)$ and so 
$$\int_{\pk2{z_j}} |e^{(j)}(z)|^qdV(z)\leqs |\rho(z_j)|^{-q\frac{i_{1,j}+i_{2,j}}2} \int_{\pk4{z_j}} |g(z)|^qdV(z).$$
In particular, for all $\alpha$ and $\beta $ and all $z\in\pk1{z_j}$, we have
\begin{eqnarray}
\left|\diffp{^{\alpha+\beta}e^{(j)} }{{z^*_1}^\alpha\partial{z^*_2}^\beta}(z) \right|&\leqs&|\rho(z_j)|^{-q\frac{i_{1,j}+i_{2,j}}2-\alpha-\frac\beta2}. \label{eq102} 
\end{eqnarray}
The inequalities (\ref{eq101}) and (\ref{eq102}) imply that the hypothesis (\ref{lp3}) of Theorem \ref{th4} is satisfied by $\hat g^{(j)}_2$ for some large $N$, the same is also true for $\hat g^{(j)}_1$.\\
Now, on ${\cal U}_2^{(j)}$, we have by Lemma \ref{lem1}~:
$$\left|\frac{P_1^{(j)}(z) e^{(j)}(z)}{Q^{(j)}_2(z)}\right| \leqs \frac1{|\rho(z_j)|} \int_{(\frac72\kappa |\rho(z_j)|)^{\frac12}\leq |\xi|\leq(4\kappa |\rho(z_j)|)^{\frac12}}
\frac{|g(\zeta^*_1,\xi)|}{\max(|f_1(\zeta^*_1,\xi )|,|f_2(\zeta_1^*,\xi)|)} dV(\xi )$$ and so 
\begin{eqnarray*}
{\int_{{\cal U}_2\cap\pk1{z_j}}\left| \frac{P_1^{(j)}(z) e^{(j)}(z)}{Q^{(j)}_2(z)}\right|^q dV(z)}
\hskip-6pt&\leqs&\hskip-6pt \int_{\pk4{z_j}}\hskip-2pt
\left(\frac{|g(\zeta^*_1,\xi)|}{\max(|f_1(\zeta^*_1,\xi )|,|f_2(\zeta_1^*,\xi)|)}\right)^q\hskip -4pt dV(\xi ) 
\end{eqnarray*}
and since $(\pk1{z_j})_{j\in\nn}$ is a $\kappa$-covering~:
\begin{eqnarray}
{\sum_{j\in\nn}\int_{{\cal U}_2\cap\pk1{z_j}}\left| \frac{P_1^{(j)}(z) e^{(j)}(z)}{Q^{(j)}_2(z)}\right|^q dV(z)}
&\leqs&\left\| \frac{g}{\max(|f_1|,|f_2|)}\right\|_{L^q(D)}^q. \label{eq103}
\end{eqnarray}
Since for all $\alpha,\beta \in\nn$, $\left|\diffp{^{\alpha+\beta}\chi^{(j)}_2}{\overline{\zeta^*_1}^\alpha\partial \overline{\zeta^*_2}^\beta} (z) \right|\leqs |\rho(z_j)|^{-\alpha-\frac\beta2}$, (\ref{eq103}) and (\ref{eq100}) imply that $(\hat g^{(j)}_2)_{j\in\nn}$ satisfy the assumption (\ref{lp2}) of Theorem \ref{th4} that we can therefore apply.\qed

\end{document}